\documentclass[12pt]{amsart}
\usepackage{amsfonts, amssymb, amsmath, amsthm, eucal, latexsym, array}
\usepackage{fullpage}
\usepackage{verbatim}

\usepackage{color}
\definecolor{darkblue}{rgb}{0,0,.5}
\definecolor{darkgreen}{rgb}{.2,0.5,.2}
\usepackage[colorlinks=true,linkcolor=darkblue,urlcolor=violet,citecolor=magenta]{hyperref}

\numberwithin{equation}{section}

\input{cyracc.def}

\newtheorem{thm}{Theorem}[section]
\newtheorem{conj}[thm]{Conjecture}
\newtheorem{lm}[thm]{Lemma}
\newtheorem{cl}[thm]{Corollary}
\newtheorem{prop}[thm]{Proposition}

\theoremstyle{remark}
\newtheorem{ex}[thm]{Example}
\newtheorem{rmk}[thm]{Remark}

\theoremstyle{definition}
\newtheorem{df}[thm]{Definition}

\newcommand{\gt}{\mathfrak}

\newcommand{\GL}{{\rm GL}}

\newcommand{\id}{{\rm id}}
\newcommand{\ind}{{\rm ind\,}}

\newcommand{\rk}{\mathrm{rk\,}}

\newcommand{\Lie}{\mathrm{Lie\,}}
\newcommand{\Ker}{{\rm Ker\,}}

\newcommand{\gr}{\mathrm{gr\,}}

\newcommand{\Ann}{\mathrm{Ann}}
\newcommand{\ad}{\mathrm{ad}}

\newcommand {\cM}{{\mathcal M}}

\newcommand {\cS}{{\mathcal S}}

\newcommand{\oA}{\overline{{\mathcal A}}}
\newcommand{\oC}{\overline{{\mathcal C}}}

\newcommand {\mK}{{\mathbb C}}

\renewcommand{\le}{\leqslant}
\renewcommand{\ge}{\geqslant}

\newfam\eusfam%
\font\euszw=eusm10 scaled 1200%
\font\eusac=eusm7 scaled 1200%
\font\eusacc=eusm7 scaled 1000%
\textfont\eusfam=\euszw\scriptfont\eusfam=\eusac%
\scriptscriptfont\eusfam=\eusacc%
\newcommand{\End}{{\rm{End}\ts}}
\newcommand{\Det}{{\rm{Det}}}
\newcommand{\Per}{{\rm{Per}}}

\newcommand{\pr}{^{\tss\prime}}
\newcommand{\non}{\nonumber}

\newcommand{\wh}{\widehat}
\newcommand{\ot}{\otimes}
\newcommand{\la}{\lambda}

\newcommand{\al}{\alpha}

\newcommand{\ga}{\gamma}
\newcommand{\Ga}{\Gamma}

\newcommand{\si}{\sigma}
\newcommand{\vp}{\varphi}
\newcommand{\de}{\delta}

\newcommand{\om}{\omega}

\newcommand{\hra}{\hookrightarrow}
\newcommand{\ve}{\varepsilon}
\newcommand{\ts}{\,}

\newcommand{\di}{\partial}

\newcommand{\tss}{\hspace{1pt}}

\newcommand{\U}{{\mathcal U}}

\newcommand{\CC}{\mathbb{C}\tss}

\newcommand{\Ac}{\mathcal{A}}

\newcommand{\Sc}{\mathcal{S}}

\newcommand{\Fc}{\mathcal{F}}

\newcommand{\Uc}{\mathcal{U}}

\newcommand{\Wc}{\mathcal{W}}

\newcommand{\gl}{\mathfrak{gl}}

\newcommand{\Pf}{{\rm Pf}}

\newcommand{\oa}{\mathfrak{o}}
\newcommand{\spa}{\mathfrak{sp}}

\newcommand{\ggot}{\mathfrak{g}}

\newcommand{\z}{\mathfrak{z}}

\newcommand{\sll}{\mathfrak{sl}}

\newcommand{\tr}{ {\rm tr}}

\newcommand{\sgn}{ {\rm sgn}\ts}

\newcommand{\Sym}{\mathfrak S}
\newcommand{\fand}{\quad\text{and}\quad}
\newcommand{\Fand}{\qquad\text{and}\qquad}

\newcommand{\bth}{\begin{thm}}
\renewcommand{\eth}{\end{thm}}
\newcommand{\bpr}{\begin{prop}}
\newcommand{\epr}{\end{prop}}
\newcommand{\ble}{\begin{lm}}
\newcommand{\ele}{\end{lm}}
\newcommand{\bco}{\begin{cl}}
\newcommand{\eco}{\end{cl}}
\newcommand{\bex}{\begin{ex}}
\newcommand{\eex}{\end{ex}}
\newcommand{\bre}{\begin{rmk}}
\newcommand{\ere}{\end{rmk}}
\newcommand{\bcj}{\begin{conj}}
\newcommand{\ecj}{\end{conj}}

\newcommand{\q}{\mathfrak{q}}

\newcommand{\bal}{\begin{aligned}}
\newcommand{\eal}{\end{aligned}}
\newcommand{\beq}{\begin{equation}}
\newcommand{\eeq}{\end{equation}}
\newcommand{\ben}{\begin{equation*}}
\newcommand{\een}{\end{equation*}}

\newcommand{\bpf}{\begin{proof}}
\newcommand{\epf}{\end{proof}}

\def\beql#1{\begin{equation}\label{#1}}

\begin{document}
\hfill {\scriptsize August 3, 2019} 
\vskip1ex

\title
{Quantisation and nilpotent limits of Mishchenko--Fomenko subalgebras}
\author[A.\,Molev]{Alexander Molev}
\address[A.\,Molev]
{School of Mathematics and Statistics,
University of Sydney,
NSW 2006, Australia}
\email{alexander.molev@sydney.edu.au}
\author[O.\,Yakimova]{Oksana Yakimova}
\address[O.\,Yakimova]{Universit\"at zu K\"oln,
Mathematisches Institut, Weyertal 86-90, 50931 K\"oln, Deutschland}
\email{yakimova.oksana@uni-koeln.de}
\begin{abstract}
For any simple Lie algebra $\ggot$ and an element $\mu\in\ggot^*$,
the corresponding commutative subalgebra $\Ac_{\mu}$ of $\U(\ggot)$ is defined
as a homomorphic image of the Feigin--Frenkel centre associated with $\ggot$.
It is known that when $\mu$ is regular this subalgebra solves Vinberg's
quantisation
problem, as the graded image of $\Ac_{\mu}$ coincides with the
Mishchenko--Fomenko subalgebra $\overline\Ac_{\mu}$ of $\Sc(\ggot)$. By a conjecture of
Feigin, Frenkel, and Toledano Laredo, this property extends to an arbitrary element $\mu$.
We give sufficient conditions on $\mu$ which imply the property.
In particular, this proves the conjecture in type {\sf C}
and gives a new proof in type {\sf A}. We show that the algebra
$\Ac_{\mu}$ is free in both cases and produce its generators in an explicit form.
Moreover, we prove that in all classical types generators of $\Ac_{\mu}$
can be obtained via the canonical symmetrisation map from certain generators of
$\overline\Ac_{\mu}$. The symmetrisation map is also used to produce free generators
of nilpotent limits of the algebras $\Ac_{\mu}$ and give a positive solution
of Vinberg's problem for these limit subalgebras.
\end{abstract}
\maketitle

\section*{Introduction}

The universal enveloping algebra $\U(\q)$ of a Lie algebra $\q$ is equipped
with a canonical filtration so that the associated graded algebra
is isomorphic to the symmetric algebra $\Sc(\q)$. The commutator on $\q$
induces the Lie--Poisson bracket on $\Sc(\q)$ defined by taking $\{X,Y\}$
to be equal to the commutator of $X,Y\in\q$ and then extending the bracket to the entire $\Sc(\q)$
by the Leibniz rule.
If $\Ac$ is a commutative subalgebra of $\U(\q)$, then its graded image $\gr\Ac$
is a Poisson-commutative subalgebra of $\Sc(\q)$. The {\em quantisation problem}
for a given Poisson-commutative subalgebra $\overline\Ac$ of $\Sc(\q)$ is to find
a commutative subalgebra $\Ac$ of $\U(\q)$ with the property $\gr\Ac=\overline\Ac$.

In the case where $\q=\ggot$ is a finite-dimensional simple Lie algebra over $\CC$,
a family of commutative subalgebras of $\U(\ggot)$ can be constructed with the use of the
associated {\em Feigin--Frenkel centre} $\z(\wh\ggot)$ which
is a commutative subalgebra
of $\U\big(t^{-1}\ggot[t^{-1}]\big)$. Given any $\mu\in\ggot^*$ and
a nonzero $z\in\CC$,
the image of $\z(\wh\ggot)$ with respect to the
homomorphism
\ben
\varrho^{}_{\ts\mu,z}\!:\U\big(t^{-1}\ggot[t^{-1}]\big)\to \U(\ggot),
\qquad X\tss t^r\mapsto X z^r+\de_{r,-1}\ts\mu(X),\quad X\in\ggot,
\een
is a commutative subalgebra $\Ac_{\mu}$ of $\U(\ggot)$ which is independent of $z$.
This subalgebra was used by Rybnikov~\cite{r:si} and
Feigin, Frenkel, and Toledano Laredo~\cite{fft:gm} to give a positive solution
of {\em Vinberg's quantisation problem} for regular $\mu$. Namely, the graded image $\gr\Ac_{\mu}$
in the symmetric algebra $\Sc(\ggot)$ turns out to coincide with the
{\em Mishchenko--Fomenko subalgebra} $\overline\Ac_{\mu}$ \cite{mf:ee} which is generated by
all $\mu$-{\em shifts} of the $\ggot$-invariants of $\Sc(\ggot)$; a precise definition
is recalled in Section~\ref{sec-MF} below. It was conjectured in \cite[Conjecture~1]{fft:gm} that
the property $\gr\Ac_{\mu}=\overline\Ac_{\mu}$ extends to all $\mu\in\ggot^*$
(it clearly holds for $\mu=0$). The conjecture was confirmed in \cite{fm:qs} for type {\sf A}.

Our first main result is a proof of the FFTL-conjecture for type {\sf C}. The same approach
can be used in type {\sf A} which leads to another proof of the conjecture.
It is known by \cite[Proposition~3.12]{fft:gm} that the inclusion
$\overline\Ac_{\mu}\subset\gr\Ac_{\mu}$ holds
for any simple Lie algebra $\ggot$ and any $\mu\in\ggot^*$.
Our argument relies on this fact and is based on a general
result establishing a maximality property
of the Mishchenko--Fomenko subalgebra $\overline\Ac_{\ga}$
associated with an arbitrary Lie algebra $\q$ and an element $\ga\in \q^*$.
In more detail, we show that under certain additional assumptions,
$\overline\Ac_{\ga}$ is a maximal Poisson-commutative
subalgebra of the algebra of $\q_{\ga}$-invariants $\Sc(\q)^{\q_{\ga}}$, where
$\q_{\ga}$ denotes the stabiliser of $\ga$
(see Theorem~\ref{thm-max}\tss(ii) below). This property is quite analogous to the main result of
\cite{codim3} establishing the maximality of the Mishchenko--Fomenko subalgebra $\overline\Ac_{\mu}$
in $\Sc(\ggot)$ for regular $\mu$.

Applying the results of \cite{ppy}, for any given $\mu\in\ggot^*$ we then produce
families of free generators of
the algebra $\Ac_{\mu}$
in types {\sf A} and {\sf C} in an explicit form.
This provides a new proof of the corresponding results of \cite{fm:qs} in type {\sf A}.

As another principal result of the paper, we show that
the free generators of $\Ac_{\mu}$ can be obtained via the canonical symmetrisation map;
see \eqref{sym} below. This map was used by Tarasov~\cite{t:cs}
to construct a commutative subalgebra of $\U(\gl_N)$ quantising the
Mishchenko--Fomenko subalgebra $\overline\Ac_{\mu}\subset \Sc(\gl_N)$.
By another result of Tarasov~\cite{t:ul}, that commutative subalgebra of $\U(\gl_N)$
coincides with $\Ac_\mu$ if $\mu$ is regular semisimple.
We extend these properties of the symmetrisation map
to all classical Lie algebras $\ggot$
by showing that the algebra of invariants $\Sc(\ggot)^{\ggot}$ admits
a family of free generators such that the images of their $\mu$-shifts
with respect to the symmetrisation map generate the algebra $\Ac_{\mu}$
for any $\mu$. The respective generators of $\Ac_{\mu}$ are given
explicitly in the form of symmetrised minors or permanents; see Theorems~\ref{thm:conja}
and \ref{thm:conjbcd}. We state as a conjecture that
free generators of $\Sc(\ggot)^{\ggot}$ with the same properties exist for
all simple Lie algebras (Conjecture~\ref{conj:sym}).

By the work of Vinberg~\cite{v:sc} and Shuvalov~\cite{Vitya},
new families
of Poisson-commu\-ta\-tive subalgebras of $\Sc(\ggot)$ of maximal transcendence degree
can be constructed by taking
certain limits of the Mishchenko--Fomenko subalgebras; see also \cite{V14}.
For instance, the graded image  of the
Gelfand--Tsetlin subalgebra ${\mathcal{GT}}(\gl_N)\subset\U(\gl_N)$ is a Poisson commutative
subalgebra of $\Sc(\gl_N)$ which does not occur as $\overline\Ac_{\mu}$
for any $\mu$. However, it can be obtained by choosing a parameter-depending family $\mu(t)$
and taking an appropriate limit of $\overline\Ac_{\mu(t)}$ as $t\to 0$.
We show that the Vinberg--Shuvalov limit subalgebras admit a quantisation.
In particular, in the case of the symplectic Lie algebra $\ggot=\spa_{2n}$
this leads to a construction of a Gelfand--Tsetlin-type subalgebra ${\mathcal{GT}}(\spa_{2n})$.
This is a maximal commutative subalgebra of $\U(\spa_{2n})$ which contains
the centres of all universal enveloping algebras $\U(\spa_{2k})$
with $k=1,\dots,n$ associated with the subalgebras
of the chain $\spa_{2}\subset\dots\subset\spa_{2n}$.
There is evidence to believe  
that this new subalgebra ${\mathcal{GT}}(\spa_{2n})$ is useful
in the representation theory of $\spa_{2n}$. 
In particular, it can be applied
to separate multiplicities in the reduction $\spa_{2n}\downarrow\spa_{2n-2}$ \cite{cris} thus leading to
a new basis for each finite-dimensional irreducible representation.

We also consider certain versions of the limit subalgebras
which are different from those of \cite{Vitya}, but arise within the general scheme
described in \cite{V14}. We
give a solution of the quantisation problem for
these {\em nilpotent limit subalgebras}
in types {\sf A} and {\sf C};
see Proposition~\ref{cor:limvin}.

It was already pointed out by Tarasov~\cite{t:ul} that the symmetrisation map
commutes with taking limits thus allowing one to quantise
the limit Poisson-commutative subalgebras of $\Sc(\ggot)$.
Therefore, the quantisations can be obtained equivalently either by
applying the symmetrisation map, or by taking the nilpotent limits of the
subalgebras $\Ac_{\mu}$.

As a consequence of the nilpotent limit construction,
we get a solution of Vinberg's quantisation problem
for centralisers of nilpotent elements in types {\sf A} and {\sf C}.
In their recent work \cite{ap}
Arakawa and Premet extended the approach of \cite{r:si} and
\cite{fft:gm} by replacing the Feigin--Frenkel centre with
the centre of the affine $\Wc$-algebra associated with a simple Lie algebra $\ggot$
and a nilpotent element $e\in\ggot$. Under certain restrictions on the data, they
produce a positive solution of Vinberg's problem for the centralisers $\ggot_e$.
It appears to be likely that their solution coincides with ours based on the nilpotent limits;
see Conjecture~\ref{conj:ap}.

Symmetric invariants of centralisers have been extensively studied at least since \cite{ppy}.
Certain polynomials ${^e\!H}\in\cS(\gt g_e)$ are defined in that paper
via the restriction to a Slodowy slice.
Notably, these elements are related to the $e$-shifts of $H$; see Lemma~\ref{l-1}.
Let $H_1,\dots,H_n\in\cS(\gt g)$ with $n=\rk\gt g$ be a set of homogeneous generating symmetric invariants.
Then $\sum\limits_{i=1}^n \deg {^e\!H}_i\le {\bf b}(\gt g_e)$,
where ${\bf b}(\gt g_e)$ is a certain integer related to $\gt g_e$.
This inequality is one of the crucial points in \cite{ppy} and it is proven 
via finite ${\mathcal W}$-algebras. We found a more direct line of argument, which works for Lie
algebras of {\em Kostant type}; see Lemma~\ref{l-2}.

Our ground field is $\mK$. However, since semisimple Lie algebras are defined over $\mathbb Z$,
it is not difficult to deduce that the main results are valid over any field of characteristic zero.

{\it Acknowledgement.} 
The work on the paper was completed during the first author's visit
to the {\em Universit\"at zu K\"oln}. He is grateful
to the {\em Mathematisches Institut}
for the warm hospitality. The second author is 
funded by the Deutsche Forschungsgemeinschaft (DFG, German Research Foundation) --- project number 330450448. We acknowledge also the support of the Australian Research
Council,
grant DP150100789.

\section{Preliminaries on Lie--Poisson structures}
Let $Q$ be a non-Abelian connected algebraic group and $\gt q=\Lie Q$ its Lie algebra.
For $\gamma\in\gt q^*$, let $\hat \gamma$ be the corresponding skew-symmetric form on $\gt q$ given by
$\hat \gamma(\xi,\eta)=\gamma\big([\xi,\eta]\big)$. Note that the kernel of $\hat\gamma$ is equal to the
stabiliser
\ben
\gt q_\gamma=\{\xi\in\gt q\mid \ad^*(\xi)\gamma=0\}.
\een

We will identify the symmetric algebra
$\Sc(\q)$ with the algebra $\CC[\q^*]$ of polynomial functions on $\q^*$.
Suppose that $\dim\gt q=r$ and choose a basis $\{\xi_1,\ldots,\xi_r\}$ of $\gt q$.
Let $\{x_1,\ldots,x_r\}$ be
the dual basis of $\gt q^*$.
Let
\beql{pi}
\pi=\sum\limits_{i<j}[\xi_i,\xi_j]\ts x_i\wedge x_j
\eeq
be the Poisson tensor (bivector) of $\gt q$. It is a global section of
$\Lambda^2 T \gt q^*$ and at each point $\gamma\in\gt q^*$ we have $\pi(\gamma)=\hat{\gamma}$.
Let $dF$ denote the differential of $F\in\cS(\gt q)$ and $d_\gamma F$ denote the differential of $F$ at $\gamma\in\gt q^*$.
A well-known property of
$\pi$ is that $$\{F_1,F_2\}(\gamma)=\pi(\gamma)(d_\gamma F_1,d_\gamma F_2)$$ for all $F_1,F_2\in\cS(\gt q)$.
As defined by Dixmier, the {\it index}
of  $\gt q$ is the number
\ben
\ind\gt q=\min_{\gamma\in\gt q^*} \dim\gt q_\gamma
=\dim\gt q-\max_{\gamma\in\gt q^*}\dim(Q\gamma)
=\dim\gt q-\max_{\gamma\in\gt q^*}\rk\pi(\gamma).
\een
For a reductive $\gt g$, one has $\ind\gt g=\rk\gt g$. In this case, $(\dim\gt g+\rk\gt g)/2$ is the dimension
of a Borel subalgebra of $\gt g$. For an arbitrary $\gt q$, set
${\bf b}(\gt q)=(\ind\gt q+\dim\gt q)/2$.
Observe that  $\ind\gt q_\gamma\ge \ind\gt q$ for all $\gamma\in\gt q^*$ by Vinberg's inequality \cite[Sect.\,1]{Dima2}.

One defines the  {\it singular set} $\gt q^*_{\rm sing}$ of $\gt q^*$ by
$$
\gt q^*_{\rm sing}=\{\gamma\in\gt q^* \mid \dim\gt q_\gamma>\ind \gt q\}.
$$
Set also $\gt q^*_{\rm reg}=\gt q^*\setminus \gt q^*_{\rm sing}$.
Further, $\gt q$ is said to have
the codim-$m$ property (or to satisfy the codim-$m$ condition), if
$\dim\gt q ^*_{\rm sing}\le \dim\gt q-m$.
A reductive Lie algebra $\gt g$ satisfies the codim-3 condition \cite{K}.
Recall that an open subset is called {\it big} if its complement does not contain divisors.
The codim-2 condition holds for $\gt q$ if and only if $\gt q^*_{\rm reg}$ is big.

Suppose that $\gamma\in\gt q^*_{\rm reg}$. Then
$$
\dim\gt q_\gamma=\ind\gt q\le \ind\gt q_\gamma.
$$
Therefore $\ind\gt q_\gamma=\dim\gt q_\gamma=\ind\gt q$ and $\gt q_\gamma$ is a commutative Lie algebra.

\subsection{Transcendence degree bounds}
\label{subsec:tdb}

For any subalgebra  $A\subset \cS(\gt q)$ and any $x\in\gt q^*$ set
\ben
d_xA=\left<d_x F \mid F\in A\right>_{\mK}\subset T^*_x\gt q^*.
\een
Then  we have
${\rm tr.deg}\,A=\max\limits_{x\in\gt q^*}\dim d_xA$.
If $A$ is Poisson-commutative, then $\hat x(d_xA,d_xA)=0$ for each $x\in\gt q^*$ and
thereby
\ben
{\rm tr.deg}\,A\le \frac{\dim\gt q-\ind\gt q}{2}+\ind\gt q={\bf b}(\gt q).
\een
Applying a result of Borho and Kraft \cite[Satz~5.7]{bk-GK}, one obtains that
${\rm tr.deg}\,{\mathcal A}\le {\bf b}(\gt q)$ for any commutative subalgebra
${\mathcal A}\subset {\mathcal U}(\gt q)$. 

For any subalgebra $\gt l\subset \gt q$, let $\cS(\gt q)^{\gt l}$ denote the
{\it Poisson centraliser} of $\gt l$, i.e.,
$$
\cS(\gt q)^{\gt l}=\{F\in\cS(\gt q) \mid \{\xi,F\}=0\quad  \text{for all}\quad \xi\in\gt l\}.
$$
If $\gt l=\Lie L$ and $L\subset Q$ is a connected subgroup, then
$\cS(\gt q)^{\gt l}$ coincides with the subalgebra of $L$-invariants
\ben
\cS(\gt q)^L=\{F\in\cS(\gt q)\mid gF=F \quad  \text{for all}\quad
g\in L\}.
\een

\begin{prop}\label{l-tr.deg}
Let $A\subset \cS(\gt q)^{\gt l}$ be a Poisson-commutative subalgebra.
Then $${\rm tr.deg}\,A\le \frac{1}{2}\ts
(\dim\gt q-\ind\gt q-\dim\gt l+\ind\gt l)+\ind\gt q={\bf b}(\gt q)-{\bf b}(\gt l)+\ind\gt l.$$
\end{prop}
\begin{proof}
For any point $x\in\gt q^*$, we have $\hat x(d_xA,d_xA)=0$ and $\hat x(\gt l,d_xA)=0$.
For a generic point $x$, the form $\hat x$ is of rank $\dim\gt q{-}\ind\gt q$ and $x|_{\gt l}\in\gt l^*_{\rm reg}$,
therefore the restriction of $\hat x$ to $\gt l$ is of rank $\dim\gt l{-}\ind\gt l$.
The quotient $d_xA/(d_xA\cap\gt q_x)$ is an isotropic subspace of $\gt q/\gt q_x$
orthogonal to $\gt l/(\gt l\cap \gt q_x)$. Knowing that the rank of $\hat x$
on $\gt l$ is equal to $\dim\gt l{-}\ind\gt l$, we can conclude that
the dimension of this quotient
is bounded by  $\frac{1}{2}\dim(Qx)-\frac{1}{2}(\dim\gt l-\ind\gt l)$.
Thereby
\ben
\dim d_xA \le \frac{1}{2}\tss\big((\dim\gt q-\ind\gt q)-(\dim\gt l-\ind\gt l)\big)
+\ind\gt q.
\een
This completes the proof.
\end{proof}

The inequality ${\rm tr.deg}\,A\le {\bf b}(\gt q)-{\bf b}(\gt l)+\ind\gt l$ has 
a different, less elementary proof. 
 According to \cite{sad}, there is a Poisson-commutative
subalgebra ${^{\gt l}\!A}\subset\cS(\gt l)$ such that  ${\rm tr.deg}\ts{^{\gt l}\!A}={\bf b}(\gt l)$.
Clearly $\{A,{^{\gt l}\!A}\}=0$ and therefore $\dim(d_xA+d_x{^{\gt l}\!A})\le{\bf b}(\gt q)$.
In addition $\dim(d_x A\cap d_x{^{\gt l}\!A})\le \ind\gt l$ for
generic $x\in\gt q^*$.

\begin{ex} {\sf (i)} Suppose that $\gt q=\gt{gl}_N$ and $\gt l=\gt{gl}_{N{-}1}$.
Then ${\rm tr.deg}\,A\le 2N-1$ for any Poisson-commutative subalgebra $A\subset\cS(\gt q)^{\gt l}$.
More generally, if $\gt l=\gt{gl}_m$ then
\ben
{\rm tr.deg}\,A\le \frac{N(N{+}1)}{2}-\frac{m(m{+}1)}{2}+m.
\een
{\sf (ii)} Take $\gt q=\gt{sp}_{2n}$ and $\gt l=\gt{sp}_{2m}$. Then ${\rm tr.deg}\,A\le n(n{+}1)-m(m{+}1)+m$.
\end{ex}

\begin{ex} Suppose that  $\gt l=\gt q_\gamma$.
If $\ind\gt q_\gamma=\ind\gt q$, then the bound of Proposition~\ref{l-tr.deg} is simpler
$${\rm tr.deg}\,A\le \frac{1}{2}(\dim\gt q-\dim\gt q_\gamma)+\ind\gt q=\frac{1}{2}\dim(Q\gamma)+\ind\gt q.$$
If $\ind\gt q_\gamma\ne\ind\gt q$, then $\ind\gt q_\gamma> \ind\gt q$ and
the bound is larger than the sum above. 
\end{ex}

Checking the equality $\ind\gt q_\gamma=\ind\gt q$ is an intricate task. It does not hold for all
Lie algebras, see e.g. \cite[Ex.~1.1]{py-gib}.

\begin{rmk}\label{comp}
The symplectic linear algebra calculation in the proof of Proposition~\ref{l-tr.deg}
brings up the following observation.
Let $A\subset \cS(\gt q)^{\gt l}$ be a Poisson-commutative subalgebra with the maximal
possible transcendence degree. Suppose that $x\in\gt q^*$ is generic, in particular   $\dim d_xA={\rm tr.deg}\,A$. Then
the orthogonal complement of $d_x A$ w.r.t. $\hat x$ equals
$\gt l+d_x A + \ker\hat x$.
\end{rmk}

\subsection{Lie algebras of Kostant type}
\label{subsec:ggs}

Throughout the paper we will use the direction derivatives $\di_{\ga}H$ of elements $H\in \Sc(\q)=\CC[\q^*]$
with respect to $\ga\in\q^*$ which are defined by
\ben
\di_{\ga}H(x)=\frac{d}{dt}\tss H(x+t\tss\ga)\Big|_{t=0}.
\een

Given a nonzero $\gamma\in\gt q^*$  we fix
a decomposition $\gt q=\mK y\tss\oplus\tss\Ker\gamma$, where $\gamma(y)=1$. For each
nonzero $H\in \cS(\gt q)$, we have a decomposition
$H=y^mH_{[m]}+y^{m{-}1}H_{[m{-}1]}+\dots +y H_{[1]}+H_{[0]}$, where $H_{[m]}\ne 0$ and
$H_{[i]}\in\cS(\Ker\gamma)$ for every $i$. Following \cite{ppy} set $^\gamma\!H=H_{[m]}$.
Note that $^\gamma\!H$ does not depend on the choice of $y$. Note also that $H_{[m]}\in\mK$
if and only if $H(\gamma)\ne 0$.

We will denote by $Q_\gamma$ the stabiliser of $\ga$ in $Q$ with respect to the
coadjoint action.

\begin{lm}\label{l-1}
Suppose that $H\in\cS(\gt q)^{\gt q}$. Let $m$ and $H_{[m]}={^\gamma\!H}$ be as above.
Then $^\gamma\!H\in\cS(\gt q_\gamma)^{Q_\gamma}$. Further\-more,
$\partial^m_\gamma H=m!H_{[m]}$ and  for all $k\ge 0$, we have $\partial^k_\gamma H=0$ if and only if $k>m$.
\end{lm}
\begin{proof}
We repeat the argument of \cite[Appendix]{ppy}. Suppose that
$H_{[m]}\not\in\cS(\gt q_\gamma)$. Then there is $\xi\in\gt q$ such that
$\{\xi,H_{[m]}\}=y\tilde H+\tilde H_0$ with $\tilde H\ne 0$ and $\tilde H, \tilde H_0\in\cS(\Ker\gamma)$. Since
$\deg _y \{\xi,y^m\} \le m$ and $\deg _y \{\xi,y^d H_{[d]}\} \le d{+}1$ for each $d$, we see that $\{\xi,H\}\ne 0$, a contradiction.
Observe that  $\Ker\gamma$ is a $Q_\gamma$-stable subspace and that
$Q_\gamma y \in y + \Ker\gamma$. Since $H$ is a $Q$-invariant and hence also a $Q_y$-invariant,
$H_{[m]}$ is a $Q_\gamma$-invariant as well.

The statements concerning derivatives follow from the facts that $\partial_\gamma y=1$ and that
$\partial_\gamma H_{[d]}=0$ for each $d$.
\end{proof}

Recall that $\xi_1,\dots,\xi_r$ is a basis of $\q$ and let
$n=\ind\gt q$.  Using notation \eqref{pi}, for any $k>0$, set
$$
\Lambda^k\pi=\,\underbrace{\pi\wedge\pi\wedge\ldots\wedge \pi}_{k \ \scriptstyle{\mathrm{factors}}}
$$
and  regard it as
an element of $\cS^k(\gt q){\otimes}\Lambda^{2k}\gt q^*$.
Then $\Lambda^{(r-n)/2}\pi\ne 0$ and all higher
exterior powers of $\pi$ are zero. We have $dF\in \cS(\gt q){\otimes}\gt q$ for each $F\in\cS(\gt q)$.
We will also regard $dF$  as a
 differential 1-form on $\gt q^*$.
Take $H_1,\ldots,H_n\in\cS(\gt q)^{\gt q}$.  Then
$dH_1\wedge\ldots\wedge dH_n\in  \cS(\gt q){\otimes}\Lambda^n \gt q$.
At the same time, $\Lambda^{(r-n)/2}\pi\in\cS(\gt q){\otimes}\Lambda^{r-n}\gt q^*$.
The volume form  $\omega=\xi_1\wedge\ldots\wedge \xi_r$
defines a non-degenerate pairing between $\Lambda^n\gt q$ and
$\Lambda^{r-n}\gt q$. If $u\in \Lambda^n\gt q$ and $v\in \Lambda^{r-n}\gt q$,
then $u\wedge v=c\,\omega$ with $c\in\mK$. We write this as
$\dfrac{u\wedge v}{\omega}=c$ and let $\dfrac{u}{\omega}$ be an element of
$(\Lambda^{r-n}\gt q)^*$ such that $\dfrac{u}{\omega}(v)=\dfrac{u\wedge v}{\omega}$.
For any ${\bf u}\in \cS(\gt q){\otimes} \Lambda^n\gt q$,
we let $\dfrac{\bf u}{\omega}$ be the corresponding element of
\ben
\cS(\gt q){\otimes}(\Lambda^{r-n}\gt q)^*\cong \cS(\gt q){\otimes}\Lambda^{r-n}\gt q^*.
\een
One says that $H_1,\ldots,H_n$ {\it satisfy the Kostant identity} if
\begin{equation}\label{K1}
\frac{dH_1\wedge\ldots\wedge dH_n}{\omega}=C\Lambda^{(r-n)/2}\pi
\end{equation}
for some nonzero constant $C$. Identity~\eqref{K1} encodes the following equivalence
$$
d_\gamma H_1\wedge\ldots\wedge d_\gamma H_n \ne 0 \  \Longleftrightarrow  \
\gamma\in\gt q^*_{\rm reg}.
$$

\begin{df}[{cf. \cite[Def.~2.2]{contr}}]\label{K2}
A Lie algebra $\gt q$ is of {\it Kostant type} if $\cS(\gt q)^Q$ is freely generated by homogenous polynomials $H_1,\ldots,H_n$ that satisfy the  Kostant identity.
\end{df}

Any reductive Lie algebra is of Kostant type \cite[Thm~9]{K}. Another easy observation is that
\ben
\sum\limits_{i=1}^{n}\deg H_i = \frac{r-n}{2}+n=\bf b(\gt q)
\een
if  homogeneous invariants $H_1,\ldots,H_n$ satisfy the Kostant identity.
If $\gt q$ is of Kostant type, then any set of algebraically independent homogeneous generators of
$\cS(\gt q)^{Q}$ satisfies the Kostant identity.

\begin{df}[{cf. \cite[Sect.~2.7]{ppy}}]\label{df-ggs}
Let $H_1,\ldots,H_n\in\cS(\gt q)^Q$ be algebraically independent homogeneous elements that
satisfy the Kostant identity~\eqref{K1}. 
If 
$\sum\limits_{i=1}^{n} \deg {^\gamma\!H_i}={\bf b}(\gt q_\gamma)$, then $\{H_i\}$
is a {\em good system}  (\tss{\em g.s.}) for $\gamma$. If in addition 
the polynomials $\{H_i\}$ generate $\cS(\gt q)^Q$, then they 
form a  {\em good generating system} (\tss{\em g.g.s.}) for $\gamma\in\gt q^*$.
\end{df}

\begin{lm}\label{l-2}
Suppose that homogeneous elements $H_1,\ldots,H_n\in\cS(\gt q)^{\gt q}$ satisfy the  Kostant identity
and $\ga\in\q^*$ is such that $\ind\gt q_\gamma=\ind\gt q$. Then
\begin{itemize}
\item[({\sf i})]
$\sum\limits_{i=1}^{n} \deg {^\gamma\!H_i}\le {\bf b}(\gt q_\gamma)$;
\item[({\sf ii})]
$\{H_i\}$ is a g.s. (for $\gamma$)
if and only if the polynomials $^\gamma\!H_i$ are algebraically independent;
\item[({\sf iii})] if $\{H_i\}$ is a g.s., then
the invariants $^\gamma\!H_i$ satisfy the Kostant identity related to $\gt q_\gamma$;
\item[({\sf iv})] if
$\{H_i\}$ is a g.s.  and
$\gt q_\gamma$ has the codim-2 property, then $\cS(\gt q_\gamma)^{Q_\gamma}=\cS(\gt q_\gamma)^{\gt q_\gamma}=\mK[\{{^\gamma\!H_i}\}]$.
\end{itemize}
\end{lm}
\begin{proof}
We may suppose that the basis elements $\xi_1,\ldots,\xi_r$ of $\q$
are chosen in such a way that
the last $\dim\gt q_\gamma$ elements $\xi_i$ form a basis of  $\gt q_\gamma$.
Whenever $\xi\in\gt q_\gamma$, we have $\gamma\big([\xi,\eta]\big)=0$ for each $\eta\in\gt q$ and hence
$\deg_y [\xi,\eta]\le 0$.  The $y$-degree of $\Lambda^{(r-n)/2}\pi$ is at most $\frac{1}{2}\dim(Q\gamma)$.
We note that
$$
r-n = r-\ind\gt q_\gamma=(r-\dim\gt q_\gamma) + (\dim\gt q_\gamma-\ind\gt q_\gamma).
$$
Choosing two complementary subspaces of $\gt q/\gt q_\gamma$ that are Lagrangian w.r.t.
$\hat\gamma$, we derive
that the highest $y$-component of $\Lambda^{(r-n)/2}\pi$ is equal (up to a nonzero scalar) to
$$
y^{(r{-}\dim\gt q_\gamma)/2}(x_1\wedge\ldots \wedge x_{r{-}\dim\gt q_\gamma})\wedge \Lambda^{(\dim\gt q_\gamma-n)/2}\pi_\gamma \ne 0,
$$
where $\pi_\gamma$ is the Poisson tensor of $\gt q_\gamma$.
The  $y$-degree of the expression on the
left hand side of \eqref{K1} is at most
$\sum\limits_{i=1}^{n} (\deg H_i- \deg {^\gamma\!H_i})$. This leads to the
 inequality
\ben
\sum_{i=1}^{n} (\deg H_i-\deg {^\gamma\!H_i}) \ge \frac{1}{2}\dim(Q\gamma)
\een
which is equivalent to
\ben
\sum_{i=1}^n \deg  {^\gamma\!H_i} \le \sum_{i=1}^{n} \deg H_i -\frac{1}{2}\dim(Q\gamma) =
\frac{1}{2}(r+n-r+\dim\q_\gamma)=
 {\bf b}(\q_\gamma).
\een
We have the equality here if and only if
$$
y^{(r-\dim\gt q_\gamma)/2}\ts \frac{d\,{^\gamma\!H_1}\wedge\ldots\wedge d\,{^\gamma\!H_n}}{\omega}
$$
is the highest $y$-component of the left hand side in \eqref{K1}.
Moreover,
this is the case if and only if the polynomials ${^\gamma\!H_i}$ are algebraically independent.
Therefore ({\sf i}) and ({\sf ii}) follow.

If $\sum\limits_{i=1}^{n} \deg {^\gamma\!H_i}= {\bf b}(\gt q_\gamma)$, then
$$
y^{(r-\dim\gt q_\gamma)/2}\ts
 \frac{d\,{^\gamma\!H_1}\wedge\ldots\wedge d\,{^\gamma\!H_n}}{\omega}
 = \tilde Cy^{(r{-}\dim\gt q_\gamma)/2}(x_1\wedge\ldots
 \wedge x_{r{-}\dim\gt q_\gamma})\wedge \Lambda^{(\dim\gt q_\gamma-n)/2}\pi_\gamma
$$
for some nonzero $\tilde C\in\mK$. Writing $\omega=(\xi_1\wedge\ldots\wedge \xi_{r{-}\dim\gt q_\gamma})\wedge  \omega_\gamma$, one  proves that
$$
 \frac{d\,{^\gamma\!H_1}\wedge\ldots\wedge d\,{^\gamma\!H_n}}{\omega_\gamma}
 = \tilde C \Lambda^{(\dim\gt q_\gamma-n)/2}\pi_\gamma.
$$
Thus ({\sf iii}) is proved as well.

The Kostant identity implies that the differentials $d\,{^\gamma\!H_i}$ are linearly independent on
$(\gt q_\gamma^*)_{\rm reg}$. If $\gt q_\gamma$ has the codim-2 property, then
$\gt q^*_{\rm reg}$ is a big open subset. Thereby the homogeneous  invariants
${^\gamma\!H_i}$ generate an algebraically closed subalgebra of $\cS(\gt q_\gamma)$, see
\cite[Thm~1.1]{ppy}. Since ${\rm tr.deg}\,\cS(\gt q_\gamma)^{Q_\gamma}={\rm tr.deg}\,\cS(\gt q_\gamma)^{\gt q_\gamma}=n$, the assertion ({\sf iv)} follows.
\end{proof}

The statements of Lemma~\ref{l-2} generalise analogous assertions proven
in \cite{ppy} for $\gt q=\gt g$ reductive and $\gamma$ nilpotent.
Parts ({\sf i}), ({\sf ii}), and ({\sf iii}) of Lemma~\ref{l-2}
constitute  \cite[Thm~2.1]{ppy}, which is proven via finite ${\mathcal W}$-algebras.
Our current approach is more direct and more general.
Examples of non-reductive Lie algebras of Kostant type can be found in 
\cite{contr}. 

Suppose that  $\ind\gt q_\gamma=\ind\gt q$ and $H_1,\ldots,H_n$ is a g.s. for $\gamma$. Then
no ${^\gamma\!H_i}$ can be a constant, see Lemma~\ref{l-2}({\sf ii}).
Therefore we must have $H_i(\gamma)=0$ for each $i$.

\subsection{Sheets and limits in reductive Lie algebras}
\label{subsec:shlim}

Suppose now that $\gt q=\gt g$ is reductive and
$Q=G$.  Choose a $G$-isomorphism  $\gt g\cong\gt g^*$.
 Recall that the {\it sheets} in $\gt g$ are the irreducible components of the
 locally closed subsets $X^{(d)}=\{\xi\in\gt g\mid \dim(G\xi)=d\}$.
At a certain further point in this paper, we will have to pass from nilpotent to arbitrary elements.
To this end, sheets in $\gt g$ and the method of associated cones developed in \cite[\S\,3]{bokr} will be used.
The {\it associated cone\/} of $\mu\in\gt g^*$ is the intersection $\overline{\mK^{^\times}\!(G\mu)}\cap \gt N$,
where $\gt N$ is the nilpotent cone.
Each irreducible component of $\overline{\mK^{^\times}\!(G\mu)}\cap \gt N$ is of dimension $\dim(G\mu)$.
Let $G\gamma$ be the dense  orbit 
 in an irreducible component of the associated cone.
The set
$\mK^{^\times}\!(G\mu)$ is irreducible, hence is  contained in a  sheet. The orbit $G\gamma$ is also contained in the same sheet.
Therefore this orbit  is unique and the associated cone is irreducible.

The following statement
should be well-known in algebraic geometry. We still give a short proof.

\begin{lm}\label{curve}
There is a morphism of algebraic varieties $\tau: \mK \to \overline{\mK^{^\times}\! (G\mu)}$
with the properties $\tau(\mK\setminus\{0\})\subset \mK^{^\times}\! (G\mu)$ and
$\tau(0)=\gamma$. Moreover, $\tau$ is given by a $1$-parameter subgroup of
$\GL(\gt g^*)$.
\end{lm}
\begin{proof}
Let $f\in\gt g$ be the image of $\gamma$ under the $G$-isomorphism $\gt g^*\cong\gt g$.
We keep the same letter $\mu$ for the image of $\mu$.
Since $\gamma$ is nilpotent,  $f$ can be included into an
$\gt{sl}_2$-triple $\{e,h,f\}\subset\gt g$.
The element  $\ad(h)$ induces a $\mathbb Z$-grading on $\gt g$, where the
components are $$\gt g_m=\{\xi\in\gt g\mid [h,\xi]=m\xi\}.$$
The centraliser $G_h$ of $h$ acts on $\gt g_{-2}$ and
$G_hf$ is a dense open  subset of $\gt g_{-2}$. Since $\gt g_{\ge 0}$ is a parabolic subalgebra of $\gt g$, we have
$G\tss(\gt g_{\ge -2})=\gt g$ and $G\tss(G_0f+\gt g_{\ge -1})$ is a dense open subset of $\gt g$,
which meets $\overline{\mK^{^\times}\!(G\mu)}$.
Hence $$G\tss(G_0f+\gt g_{\ge -1})\cap \overline{\mK^{^\times}\!(G\mu)}$$ is a nonempty open
subset of $\overline{\mK^{^\times}\!(G\mu)}$ and there is $cg\mu\in G(G_0f+\gt g_{\ge -1})$ with
$g\in G$, $c\in\mK^{^\times}\!$.
We may assume that $\mu\in f+\gt g_{\ge -1}$.

Let $\{\chi(t)\mid t\in\mK^{^\times}\!\}\subset {\rm GL}(\gt g)$ be the $1$-parameter subgroup defined by
$$\chi(t)=t^2\id_{\gt g}\exp(\ad(th)).$$
Then $\lim\limits_{t\to 0}\chi(t)\mu=f$ and
$\tau$ is given by $\tau(t)=\chi(t)\mu$ for $t\ne 0$ and $\tau(0)=f$.
\end{proof}

The existence of $\tau$ and the fact that
$\dim (G\gamma)=\dim (G\mu)$  imply that
$
\text{$\displaystyle\lim_{t\to 0}\gt g_{\tau(t)}=\gt g_\gamma$,} 
$
where the limit is taken in a suitable Grassmannian.

Let $\mu=x{+}y$ be the Jordan decomposition of $\mu$ in $\gt g$,
where $x$ is semisimple and $y$ is nilpotent.
Then $\gt l=\gt g_x$ is a Levi subalgebra of $\gt g$. Set $L=\exp(\gt l)$.
By  \cite[Sect.~3]{Bor}, the nilpotent orbit $G\gamma$ in a sheet of $\mu$ is
{\it induced} from $Ly\subset\gt l$. In  standard notation,
$G\gamma={\rm Ind}_{\gt l}^{\gt g}(Ly)$. For the classical Lie algebras, the description of this
induction is given in \cite{ke83}.

\subsubsection*{Type {\sf A}}
Let $\gt g=\gl_N$.
Suppose that the distinct eigenvalues of $\mu\in\gl_N$ are $\la_1,\dots,\la_r$
and that the Jordan canonical form of $\mu$ is the direct sum of the respective Jordan blocks
\beql{jord-s}
J_{\al^{(1)}_1}(\la_1),\dots,J_{\al^{(1)}_{m(1)}}(\la_1),\dots,
J_{\al^{(r)}_1}(\la_r),\dots,J_{\al^{(r)}_{m(r)}}(\la_r),
\eeq
of sizes
$\al^{(i)}_1\geqslant \al^{(i)}_2\geqslant\dots  \geqslant \al^{(i)}_{m(i)}\geqslant 1$
for $i=1,\dots,r$.
We let $\al^{(i)}$ denote the corresponding Young diagram whose $j$-th row
is $\al^{(i)}_j$.
The nilpotent orbit $G\gamma$ 
corresponds to the partition $\Pi$
such that the row $k$ of its Young diagram is the sum of the $k$-th rows of all
diagrams $\al^{(i)}$. That is,
\beql{jord-n}
\Pi=\alpha^{(1)}+\dots+\alpha^{(r)}=\Big(\sum_{i=1}^r \al^{(i)}_1, \sum_{i=1}^r \al^{(i)}_2,\dots \Big),
\eeq
where we assume that $\al^{(i)}_j=0$ for $j>m(i)$; see \cite[\S1]{ke83} and in particular Cor.~2 there.

\subsubsection*{Type {\sf C}}
Suppose now that $\gt g=\spa_{2n}$.  Keep notation \eqref{jord-s} for the Jordan blocks of
$\mu\in\spa_{2n}$. To separate the zero eigenvalue we will
assume that $\lambda_1=0$. The case where zero is not an eigenvalue of $\mu$ will be taken care of
by the zero multiplicity $m(1)=0$. In type {\sf C}, the following additional conditions
on the parameters of the Jordan canonical form must be satisfied, cf.
\cite[\S2]{ke83} and \cite[Sect.~1]{Jan},
\begin{itemize}
\item[$\diamond$] any row of odd length in $\al^{(1)}$ must occur an even number of times;
\item[$\diamond$] for each $k>1$ there is $k'$ such that $\lambda_k=-\lambda_{k'}$
and $\al^{(k)}=\al^{(k')}$.
\end{itemize}

Define the partition $\Pi$ by \eqref{jord-n}. It may happen that
some rows of odd length in $\Pi$ occur odd number of times. We will modify $\Pi$ in order
to produce a new partition $\Pi_{\ga}$, which
corresponds to a nilpotent orbit in $\gt{sp}_{2n}$, by the  sequence of steps described in the proof of  \cite[Lemma~2.2]{ke83}. Working
from the top of the current Young diagram, consider the first row $\beta$ of odd length
which occurs an odd number of times. Remove one box from the last occurrence
of $\beta$ and add this box to the next row.
This operation is possible because the length of the next row is necessarily odd due to the
conditions on the Jordan canonical form in type {\sf C}.
Repeating the procedure will yield
a diagram $\Pi_{\ga}$ with the property that any row of odd length
occurs an even number of times. The orbit $G\gamma$ is given by the partition $\Pi_\gamma$
\cite[Prop.~3.2]{ke83}.


\section{Mishchenko--Fomenko subalgebras and their limits}\label{sec-MF}
 Take $\gamma\in\gt q^*$ and let $\overline{{\mathcal A}}_\gamma$ denote the
 corresponding
{\em Mishchenko--Fomenko subalgebra} of $\Sc(\q)$ which is generated by all
$\ga$-{\em shifts} $\di^{\tss k}_\gamma H$ with $k\ge 0$ of all elements $H\in\cS(\gt q)^{\gt q}$.
Note that $\di^{\tss k}_\gamma H$ is a constant for $k=\deg H$.
The  generators $\partial^k_\gamma H$ of
$\oA_\gamma$
can be equivalently defined by shifting
the argument of the invariants $H$.
Set
$H_{\gamma,t}(x)=H(x+t\gamma)$. Suppose that $\deg H=m$. Then
$H_{\gamma,t}(x)$ expands as 
\beql{polte}
H(x+t\gamma)
=H^{}_{(0)}(x)+ t H^{}_{(1)} (x) +  \dots+ t^m H^{}_{(m)} (x),
\eeq
where   $\ts H^{}_{(k)} = \frac{1}{k!}\partial^k_\gamma H$.
The subalgebra $\overline\Ac_{\ga}$
is generated by all elements $H^{}_{(k)}$
associated with all $\gt q$-invariants $H\in \Sc(\gt q)^{\gt q}$.
When $H$ is a homogeneous polynomial,
we will also use an equivalent form of \eqref{polte},
where we make
the substitution $x\mapsto \ga +\ts z^{-1} x$
for a variable $z$,
and expand as 
\ben
H\big(\gamma+z^{-1}x)
=z^{-m} H^{}_{(0)} (x) +\dots+z^{-1} H^{}_{(m-1)}(x) +H^{}_{(m)}(x).
\een

Take now $\lambda\in\mK$ and consider $H_{\gamma,\lambda}\in \cS(\gt q)$ with
$H_{\gamma,\lambda}(x)=H(x+\lambda\gamma)$.
A standard argument with the Vandermonde determinant  shows that
$\oA_\gamma$ is generated by $H_{\gamma,\lambda}$ with $H\in\cS(\gt q)^{\gt q}$ and $\lambda\in\mK$.
One readily sees that
\ben
\oA_\gamma \subset \cS(\gt q)^{Q_\gamma}.
\een
Clearly it suffices   to take only homogeneous $H$ and only nonzero $\lambda$ in order to generate
$\oA_\gamma$. Under these assumptions on $H$ and $\lambda$, we have
\begin{equation}\label{d-lambda}
d_\mu H_{\gamma,\lambda}=d_{\mu{+}\lambda\gamma}H=\lambda^{m{-}1} d_{\lambda^{-1}\mu{+}\gamma}H=\lambda^{m{-}1} d_{\gamma} H_{\mu,\lambda^{-1}},
\end{equation}
where $m=\deg H$.  Hence
\begin{equation}\label{dif-A}
d_\mu \oA_\gamma=d_\gamma \oA_\mu
\end{equation}
for any $\gamma,\mu\in\q^*$, where both sides
are regarded as subspaces of $\gt q$.  

\subsection{Maximal Poisson-commutative subalgebras}
\label{subsec:mpc}

 Our next goal is to formulate certain conditions that assure  that $\oA_\gamma$ is a maximal Poisson-commutative subalgebra of
$\cS(\q)^{\gt q_\gamma}$.
From now on assume that  $\gt q$ is of Kostant type with $\ind\q=n$
and that $H_1,\ldots,H_n$
are homogeneous generators of $\cS(\gt q)^{\gt q}$. 
Then the Mishchenko--Fomenko subalgebra $\oA_\gamma$
is generated by the
$\ga$-shifts $\di^{\tss k}_\gamma H_i$ with
$1\le i\le n$ and $0\le k\le \deg H_i-1$.

\begin{lm}\label{lm-bol}
Suppose that $\gt q$ is of Kostant type and has the codim-2 property. Suppose further  that
$\ind\gt q_\gamma=n$. Then
\begin{itemize}
\item[({\sf i})] ${\rm tr.deg}\,\oA_\gamma=\frac{1}{2}\dim(Q\gamma)+n$;
\item[({\sf ii})] $\dim d_x \oA_\gamma=\frac{1}{2}\dim(Q\gamma)+n$\ \  if\ \
$x+\mK\gamma\subset \gt q^*_{\rm reg}$\ \  and\ \  $x|_{\gt q_\gamma}\in(\gt q_\gamma^*)_{\rm reg}$.
\end{itemize}
\end{lm}
\begin{proof}
By \eqref{dif-A}, $d_x \oA_\gamma=d_\gamma \oA_x$.
According to  \cite[Thm~3.2 and its proof]{bol}, see also \cite[Thm~2.1]{bol},  if $x+\mK\gamma\subset \gt q^*_{\rm reg}$ and $x|_{\gt q_\gamma}\in(\gt q_\gamma^*)_{\rm reg}$, then
$d_\gamma \oA_x$ contains a subspace $U$ of dimension $\frac{1}{2}\dim(Q\gamma)$ such that
$U\cap \gt q_\gamma=0$. A part of \cite[Thm~3.2]{bol} asserts that such an element $x$ exists.
For each $H\in\cS(\gt q)^{\gt q}$ and each $\mu\in\gt q^*$, we have $d_\mu H\in\Ker\hat\mu=\gt q_\mu$.
Since $\gt q$ is of Kostant type, the differentials $d_{ux{+}\gamma} H_i$ are linearly independent
for each $u\in\mK^{^\times}\!$ and therefore their linear span is equal to $\gt q_{ux{+}\gamma}$.
In view of \eqref{d-lambda}, $\gt q_{ux{+}\gamma}\subset d_\gamma \oA_x$ and hence
$U_\gamma=\lim\limits_{u\to 0} \gt q_{ux{+}\gamma}$ is a subspace of $d_\gamma \oA_x$ as well. Clearly
$\dim U_\gamma=n$ and $U_\gamma\subset \gt q_\gamma$. It follows that
\ben
\dim d_\gamma\oA_x\ge \dim U+\dim U_\gamma=\frac{1}{2}\dim(Q\gamma)+n.
\een
Since $d_x \oA_\gamma \le \frac{1}{2}\dim(Q\gamma)+n$ by Proposition~\ref{l-tr.deg}, part ({\sf ii}) follows.
To prove ({\sf i})
recall  that
${\rm tr.deg}\,\oA_\gamma=\max_{x\in\gt q^*}\dim d_x \oA_\gamma=\frac{1}{2}\dim(Q\gamma)+\ind\gt q$.
\end{proof}



\begin{prop}\label{Prop-Bol}
Suppose that $\gt q$ satisfies the codim-2 condition and that $H_1,\ldots, H_n$ is a g.g.s. for $\gamma\in\gt q^*$.
If in addition $\ind\gt q_\gamma=n$, then $\oA_\gamma$
is freely generated by $\partial_\gamma^k H_i$ with $1\le i\le n$ and
$0\le k \le \deg H_i-\deg {^\gamma\!H_i}$.
\end{prop}

\bpf
Since
$\sum\limits_{i=1}^{n}\deg {^\gamma\!H_i}={\bf b}(\gt q_\gamma)$, there are only
$$
\sum\limits_{i=1}^{n}(\deg H_i-\deg {^\gamma\!H_i}+1)=\frac{\dim\gt q+n}{2}-\frac{\dim\gt q_\gamma+n}{2}+n=\frac{1}{2}\dim(Q\gamma)+n
$$
nonzero elements $\di_\gamma^k H_i$ and they  have to be algebraically independent by Lemma~\ref{lm-bol}({\sf i)}.
\epf

\begin{thm}\label{thm-max}
Suppose that $\gt q$  satisfies the codim-3 condition and that $H_1,\ldots,H_n$ is a g.g.s. for $\gamma$. Suppose further that
$\gt q_\gamma$ satisfies the codim-2 condition and $\ind\gt q_\gamma=n$.
Let $F_1,\ldots,F_s$ with $s=\frac{1}{2}\dim(Q\gamma)+n$ be algebraically independent homogeneous generators of
$\oA_\gamma$ such that  each $F_j$ is a $\gamma$-shift $\partial^k_\gamma H_i$ as in Proposition~\ref{Prop-Bol}.
Then
\begin{itemize}
\item[({\sf i})] the differentials of $F_j$ are linearly independent on a big open subset;
\item[({\sf ii})]   $\oA_\gamma$ is a maximal (w.r.t. inclusion) Poisson-commutative subalgebra of
 $\cS(\gt q)^{\gt q_\gamma}$.
\end{itemize}
\end{thm}
\begin{proof} ({\sf i}) We have $d_x \oA_\gamma=\left<d_xF_j \mid 1\le j\le s\right>_{\mK}$.
Whenever $\dim d_x\oA_\gamma=s$, the differentials $d_xF_j$ are linearly independent.
According to Lemma~\ref{lm-bol},
the equality  $\dim d_x \oA_\gamma=s$ holds if
$x+\mK \gamma\subset \gt q^*_{\rm reg}$ and
$x|_{\gt q_\gamma}\in (\gt q_\gamma^*)_{\rm reg}$. 
Choosing any complement of $\gt q_\gamma$ in $\gt q$, we can embed $\gt q_\gamma^*$ into
$\gt q^*$.
If $x|_{\gt q_\gamma}$ is non-regular, then $x$ belongs
to $(\gt q_\gamma^*)_{\rm sing}+\Ann(\gt q_\gamma)$.  This is a closed subset of codimension at least $2$, since $\gt q_\gamma$  satisfies the codim-2 condition.

Let us examine the property $x+\mK \gamma\subset \gt q^*_{\rm reg}$.
The desired condition on $x$ holds if
$x\in\gt q^*_{\rm reg}$ and $\gamma+cx\in\gt q^*_{\rm reg}$ for all $c\in\mK^{^\times}\!$.
The first  restriction is inessential. In order to deal with the second one,
$\gamma+\mK^{^\times}\! x\subset\gt q^*_{\rm reg}$,
we
choose $\gamma$ as the origin in $\gt q^*$ and consider the corresponding
map $$\psi\!: \gt q^*\setminus\{\gamma\} \to \mathbb P\gt q^*$$ with $\psi(x)=\mK(x{-}\gamma)$.
We have
\ben
\dim\overline{\psi(\gt q^*_{\rm sing})}\le \dim\gt q^*_{\rm sing}\le \dim\gt q-3.
\een
Hence the preimage $\psi^{-1}(\overline{\psi(\gt q^*_{\rm sing})})$ is a closed
subset of $\gt q^*\setminus\{\gamma\}$ of codimension at least $2$. Assume that $x\ne 0$. Note that
$$\gamma+\mK^{^\times}\! x = \psi^{-1}\big(\psi(x+\gamma)\big).$$
If $(\gamma+\mK^{^\times}\! x)\cap\gt q^*_{\rm sing}\ne \varnothing$,
then $\psi(x+\gamma)\in\psi(\gt q^*_{\rm sing})$ and
$x+\gamma\in \psi^{-1}\big(\psi(\gt q^*_{\rm sing})\big)$.

Since
\ben
\psi^{-1}\big(\tss\overline{\psi(\gt q^*_{\rm sing}})\big)\cup\{\gamma\}
\een
is a closed subset of $\gt q^*$ of dimension at most $\dim\gt q{-}2$, part  ({\sf i}) follows.

({\sf ii}) We have
${\rm tr.deg}\ts\oA_\gamma=\frac{1}{2}\dim(Q\gamma)+n$.
Assume on the contrary that $\oA_\gamma$ is not maximal. Then
$\oA_\gamma\subsetneq A\subset\cS(\gt q)^{\gt q_\gamma}$, where $A$ is a Poisson-commutative
subalgebra. 
In view of Proposition~\ref{l-tr.deg}, 
${\rm tr.deg}\,A\le {\rm tr.deg}\,\oA_\gamma$ and hence
$A$ is a non-trivial algebraic extension of $\oA_\gamma$.  Since  ({\sf i}) holds,
$\oA_\gamma$ is an algebraically closed subalgebra of $\cS(\gt q)$ by \cite[Thm~1.1]{ppy}.
This contradiction completes the proof.
\end{proof}


Now suppose that $\gt g=\Lie G$ is reductive. Then $\gt g$ has the codim-3 property \cite{K}.
It will be convenient to consider elements of $\gt g$ as linear functions on
$\gt g$. We have $\ind\gt g_\gamma=\rk\gt g=n$ for each $\gamma\in\gt g^*$ \cite{fan,graaf,CM}.
A nilpotent element $\gamma$ may or may not posses a good generating system \cite{ppy}.  But if an element is not nilpotent, then there is no g.g.s. for it.
Our next step will be to develop a transition from
nilpotent to arbitrary elements of $\ggot$.

We point out a few properties of directional derivatives to be used below. We have
\ben
\partial_{g\mu}(gF)=g(\partial_\mu F)
\een
for all $g\in G$, $\mu\in\gt g^*$, and $F\in\cS(\gt g)$. Hence,
$\partial_{g\mu}^k H=g(\partial_\mu^k H)$ for $H\in\cS(\gt g)^G$.
Moreover, $\partial_{t\mu}F=t\partial_\mu F$ for each $t\in\mK^{^\times}\!$.

Consider an arbitrary element $\mu\in\ggot^*$ and the associated nilpotent orbit
$G\ga\subset \ggot^*$ as defined in the first paragraph  of Section~\ref{subsec:shlim}.

\begin{prop}\label{prop-cone}
Suppose that $\gamma$ has a g.g.s. and that $\gt g_\gamma$ satisfies the codim-2 condition.
Then $\oA_\mu$ is a maximal Poisson-commutative subalgebra of $\cS(\gt g)^{\gt g_\mu}$.
\end{prop}
\begin{proof}
Let
$F_1,\ldots,F_s$ with $s=\frac{1}{2}\dim(G\gamma)+\rk\gt g=\frac{1}{2}\dim(G\mu)+\rk\gt g$
be algebraically independent homogeneous generators of $\oA_\gamma$. As in Theorem~\ref{thm-max},
we have $F_j=\partial_\gamma^k H_i$, where $0\le k \le \deg H_i-\deg{^\gamma\!H_i}$.
Set accordingly $\hat F_j=\partial_\mu^k H_i$.
By Theorem~\ref{thm-max}({\sf i}), the differentials 
$dF_j$ are linearly independent on a big open subset.
Assume on the contrary that this is not the case for the differentials of $\hat F_j$. Then
$$
\bigwedge\limits_{j=1}^s d\hat F_j = {\bf F} R,
$$
where ${\bf F}$ is a non-constant homogenous polynomial and
$R\in \cS(\gt g)\otimes \Lambda^s \gt q$ is a regular differential $s$-form that
is nonzero on a big open subset of $\gt q^*$.

Let $\tau$ be the map of  Lemma~\ref{curve}, which is constructed as an orbit of
$$\{\chi(t)=t^2\id_{\gt g^*}\exp(\ad^*(th))\mid t\in\mK^{^\times}\!\}\subset\GL(\gt g^*).$$ 
Then $\lim\limits_{t\to 0}\partial_{\tau(t)}^k H_i=\partial_{\gamma}^k H_i$ for
all $i$ and all $k$. The appearing partial derivatives can  be expressed via the group action:
 $$\partial_{\tau(t)}^k H_i=\partial_{\chi(t)\mu}^k H_i=t^{2k}\exp(\ad^*(th))(\partial_\mu^k H_i).$$
Letting $G$ act on the differential forms as well, we obtain that
$$
\bigwedge\limits_{j=1}^s d F_j = \lim\limits_{t\to 0} t^K(\exp(\ad^*(th)){\bf F}) (\exp(\ad^*(th))R),
$$
where $K\in\mathbb Z_{\ge 0}$.
Since ${\bf F}\not\in\mK$ is a homogeneous polynomial, the lowest $t$-component of
$\exp(\ad^*(th)){\bf F}$ is a non-constant homogeneous polynomial as well.
This component divides $dF_1\wedge\ldots\wedge dF_s$, a contradiction.
Thus, the differentials $d\hat F_j$ are linearly independent on a big open subset.

By \cite[Thm~1.1]{ppy}, the homogeneous elements $\hat F_1,\ldots, \hat F_s$
generate an algebraically closed subalgebra of $\cS(\gt g)$,
which is contained in $\oA_\mu$ and is of transcendence degree $s$.
Since ${\rm tr.deg}\,\oA_{\mu}=s$, these elements actually generate $\oA_{\mu}$.
By Proposition~\ref{l-tr.deg}, a Poisson-commutative extension of $\oA_{\mu}$ in
$\cS(\gt g)^{\gt g_{\mu}}$ must be algebraic and is therefore trivial.
\end{proof}

If $\gt q=\gt g$ is of type {\sf A} or {\sf C} and $\gamma$ is nilpotent, then
$\gt g_\gamma$ has the codim-2 property and there is a g.g.s. for $\gamma$, see \cite{ppy}.
Therefore, the assumptions of Theorem~\ref{thm-max} are satisfied
and we have the following.

\begin{cl}\label{A-C}
Suppose that $\gt g$ is of type {\sf A} or {\sf C}. Then
$\oA_\mu$ is a maximal Poisson-commutative subalgebra of $\cS(\gt g)^{\gt g_\mu}$ for each
$\mu\in\gt g^*$.
\qed
\end{cl}

If  $\gamma\in\gt q^*_{\rm sing}$, then ${\rm tr.deg}\,\oA_\gamma < {\bf b}(\gt q)$ \cite[Thm~2.1]{bol}, i.e., this algebra does not
have the maximal possible 
transcendence degree.  On the one hand,
this property can be a disadvantage for some applications, while on the
other, $\oA_\gamma$ Poisson-commutes with $\gt q_\gamma$ and can be extended, in many ways, to a Poisson-commutative
subalgebra of the maximal possible transcendence degree. Below we present a construction of such an extension.

\subsection{Vinberg's limits in the nilpotent case}
\label{sec-limits}

Take $\gamma\in \gt q^*_{\rm sing}$,
 $\mu\in\gt q^*_{\rm reg}$, and $u\in\mK$ and consider
the Mishchenko--Fomenko subalgebra $\oA_{\gamma{+}u\mu}$ of $\Sc(\q)$.
For $F\in \cS^N(\gt  q)$, we have $dF\in \cS^{N{-}1}(\gt q)\otimes \gt q$
and $\partial_x F=dF(\,.\,,x)$ for any $x\in \gt q^*$. Hence $\partial_{\gamma+u\mu}F=\partial_\gamma F+u\partial_\mu F$. More generally, $\partial^k_{\gamma+u\mu}F\in \cS(\gt q)[u]$.
We have
\begin{equation}\label{par}
\partial^k_{\gamma{+}u\mu}
=\partial_\gamma^k+k u\partial_\mu\partial_\gamma^{k{-}1}+\dots
+ \binom{k}{s} u^s \partial_\mu^s \partial^{k{-}s}_\gamma+\dots+u^k\partial^k_\mu.
\end{equation}
Obviously
$\lim\limits_{u\to 0} \partial^k_{\gamma+u\mu}F =\partial_\gamma^k F$.
If $\partial_\gamma^k F=0$, but $\partial_{\gamma+u\mu}^k F\ne 0$, then 
the limit
$\lim\limits_{u\to 0}  \mK\partial^k_{\gamma+u\mu}F$ still makes sense
as an element of the projective space
$\mathbb P\tss\cS(\gt q)$. This limit line is spanned
by the lowest $u$-component of $ \partial^k_{\gamma+u\mu}F$.
In the same projective sense set
$$
\oC_{\gamma,\mu}=\lim\limits_{u\to 0} \oA_{\gamma{+}u\mu}.
$$
Formally speaking, $\oC_{\gamma,\mu}$ is a subspace of $\cS(\gt q)$ such that
$$
\oC_{\gamma,\mu}\cap \cS^m(\gt q)=\lim_{u\to 0}\ts \big(\tss\oA_{\gamma{+}u\mu}\cap \cS^m(\gt q)\big)
$$
for each $m\ge 0$. In other words, 
\begin{equation}\label{V-limit-def}
\oC_{\gamma,\mu}=\left<\tss \text{lowest $u$-component of }
\boldsymbol{F}=F_0+uF_1+\dots + u^k F_k \mid k\ge 0, F_i\in\oA_{\gamma{+}u\mu} \right>_{\mK}.
\end{equation}
We call $\oC_{\gamma,\mu}$ {\it Vinberg's limit} at $\gamma$ along $\mu$, see \cite{V14}.
Note that $\oC_{\gamma,\mu}$ is a subalgebra of $\cS(\gt q)$ and that it  does depend on $\mu$.

Clearly $\oA_\gamma\subset \oC_{\gamma,\mu}$. By
 \cite[Satz~4.5]{bk-GK}, we have
${\rm tr.deg}\,\oC_{\gamma,\mu}={\rm tr.deg}\,\oA_{\gamma+u\mu}$ with a generic $u\in\mK$. Therefore
 ${\rm tr.deg}\,\oC_{\gamma,\mu}={\bf b}(\gt q)$ assuming that $\gt q$ satisfies the codim-2 condition and
has enough symmetric invariants \cite{bol}. Set $\bar\mu=\mu|_{\gt q_\gamma}$.

\begin{thm} \label{VL-nilp}
Suppose that
$\gt q$  satisfies the codim-2 condition, $\ind\gt q_\gamma=n$,
there is a g.g.s. for $\gamma$, $\gt q_\gamma$ has the codim-2 property,   and $\bar\mu\in(\gt q_\gamma^*)_{\rm reg}$.
Then $\oC_{\gamma,\mu}$ is a free algebra generated by $\oA_\gamma$ and the Mishchenko--Fomenko subalgebra  $\oA_{\bar\mu}\subset \cS(\gt q_\gamma)$.
\end{thm}
\begin{proof}
Let $H_1,\ldots,H_n$ be a g.g.s. for $\gamma$. Set $F_{i,k}=\partial_\gamma^k H_i$.
For $k= \deg H_i{-}\deg{^\gamma\!H_i}$, we have $\partial_\gamma^k H_i=k!{^\gamma\!H_i}\in\cS(\gt q_\gamma)$ and $\partial_\gamma^P H_i=0$ if $P>k$, see
Lemma~\ref{l-1}. 
Together with \eqref{par} this gives
\ben
\lim_{u\to 0} \mK\partial_{\gamma{+}u\mu}^k H_i
=\begin{cases}
\mK F_{i,k}\qquad &\text{if }\ \  k\le \deg H_i-\deg{^\gamma\!H_i},\\
\mK\partial_{\bar{\mu}}^{\bar k} ({^\gamma\!H_i}) \qquad &\text{if }\ \  k > \deg H_i-\deg{^\gamma\!H_i},
\end{cases}
\een
where $\bar k=k-(\deg H_i-\deg{^\gamma\!H_i})$. According to Proposition~\ref{Prop-Bol}, the elements
$F_{i,k}$ with $k$ satisfying the conditions
$k\le \deg H_i-\deg{^\gamma\!H_i}$ freely generate $\oA_\gamma$.

Since $H_1,\ldots,H_n$ is a g.g.s. for $\gamma$,
the elements ${^\gamma\!H_1},\ldots,{^\gamma\!H_n}$ are algebraically independent, they
fulfil the Kostant identity~\eqref{K1} and  freely
generate $\cS(\gt q_\gamma)^{\gt q_\gamma}$, see  Lemma~\ref{l-2}.
In particular,  $\gt q_\gamma$ is of Kostant type.
Applying Proposition~\ref{Prop-Bol} to $\cS(\gt q_\gamma)$, we see that
$\oA_{\bar\mu}$ is freely generated by the elements $\partial_{\bar{\mu}}^{\bar k} ({^\gamma\!H_i})$
with $0\le \bar k < \deg {^\gamma\!H_i}$.

It remains to show that there are no algebraic relations among
$F_{i,k}$ with $k\le \deg H_i{-}\deg{^\gamma\!H_i}$ and $\partial_{\bar{\mu}}^{\bar k} ({^\gamma\!H_i})$
with $1\le \bar k < \deg {^\gamma\!H_i}$. Once this is done, we will know that the
lowest $u$-components
of the generators $\partial_{\gamma{+}u\mu}^k H_i$
are algebraically independent and therefore the lowest $u$-component  of a polynomial in $u^s\partial_{\gamma{+}u\mu}^k H_i$
is a polynomial in the lowest $u$-components of $\partial_{\gamma{+}u\mu}^k H_i$.

Take  $x\in\gt q^*$. Then $U_1=\left<d_x F_{i,k}\right>_{\mK}$ is a subspace
of $\gt q$ and $\hat x(U_1,\gt q_\gamma)=0$. At the same time
$$U_2=\left<d_x \partial_{\bar{\mu}}^{\bar k} ({^\gamma\!H_i})
\mid 0\le \bar k < \deg {^\gamma\!H_i}\right>_{\mK}$$ is a subspace of $\gt q_\gamma$.  Therefore $\hat x(U_1\cap U_2,\gt q_\gamma)=0$ and
$U_1\cap U_2\subset (\gt q_\gamma)_{\bar x}$ for $\bar x=x|_{\gt q_\gamma}$.
Suppose that $x$ is a generic point, $x\in\gt q^*_{\rm reg}$ and $\bar x\in(\gt q_\gamma^*)_{\rm reg}$.
Since $\ind\gt q_\gamma=n$, we have $\dim(\gt q_\gamma)_{\bar x}=n$. Hence $\dim(U_1\cap U_2)=n$.
It follows that $U_1+U_2=U_1\oplus \tilde U_2$, where
\ben
\tilde U_2=\Big<d_x \partial_{\bar{\mu}}^{\bar k} ({^\gamma\!H_i})
\mid 1\le \bar k < \deg {^\gamma\!H_i}\Big>_{\mK}.
\een
This completes the proof.
\end{proof}

The assumptions of Theorem~\ref{VL-nilp} are satisfied in types {\sf A} and {\sf C} for all nilpotent elements $\gamma$
and for generic elements $\mu$. 

\begin{rmk}\label{A-l-max}
Suppose that $\gt g$ is of type {\sf A}. Then
$\gt g_\gamma$ has the codim-3 property \cite{surp} and thereby
$\oA_\nu$ is a maximal Poisson-commutative subalgebra of
$\cS(\gt g_\gamma)$ for each $\nu\in(\gt g_\gamma^*)_{\rm reg}$, see
\cite{codim3} and also \cite{ap}. Hence
$\oC_{\gamma,\mu}$
is a maximal Poisson-commutative subalgebra of $\cS(\gt g)$ for any nilpotent $\gamma$ and any generic $\mu$. 
\end{rmk}

\section{Quantisation and symmetrisation}\label{sec:qs}

As we recalled in the Introduction, Vinberg's quantisation problem \cite{v:sc} concerns the existence
and construction of a commutative subalgebra $\Ac_{\mu}$ of $\U(\ggot)$
with the property $\gr\Ac_{\mu}=\overline\Ac_{\mu}$.
In the case where $\mu\in\ggot^*$ is regular semisimple,
explicit constructions of the subalgebras $\Ac_{\mu}$
in the classical types were given by Nazarov and Olshanski~\cite{no:bs}
with the use of the Yangian for $\gl_N$ and the twisted Yangians
associated with the orthogonal and symplectic Lie algebras.
Positive solutions in the general case
were given by Rybnikov~\cite{r:si} for regular semisimple $\mu$
and Feigin, Frenkel and Toledano Laredo~\cite{fft:gm} for any regular $\mu$.
The solutions rely on the properties of a commutative subalgebra $\z(\wh\ggot)$ of the
universal enveloping algebra $\U\big(t^{-1}\ggot[t^{-1}]\big)$.
This subalgebra is known as
the {\em Feigin--Frenkel centre} and is defined as the centre of the universal affine vertex algebra
associated with the affine Kac--Moody algebra $\wh\ggot$
at the critical level. In particular,
each element of $\z(\wh\ggot)$ is annihilated by the
adjoint action of $\ggot$.
Furthermore,
the subalgebra $\z(\wh\ggot)$ is invariant with respect to
the derivation $T=-\di_t$ of the algebra $\U\big(t^{-1}\ggot[t^{-1}]\big)$.
By a theorem of Feigin and Frenkel~\cite{ff:ak} (see also \cite{f:lc}),
there exists a family of elements $S_1,\dots,S_n\in \z(\wh\ggot)$
(a {\em complete set of
Segal--Sugawara vectors}), where $n=\text{\rm rk}\ts\ggot$,
such that
\ben
\z(\wh\ggot)=\CC[T^{\tss r}S_l\ |\ l=1,\dots,n,\ \ r\geqslant 0].
\een
One can assume that each $S_l$ is homogeneous with respect to
the gradation on $\U\big(t^{-1}\ggot[t^{-1}]\big)$ defined by $\deg X[r]=-r$,
where we set $X[r]=X t^{\tss r}$.

For any $\mu\in\ggot^*$ and
a nonzero $z\in\CC$, the mapping
\beql{evalr}
\varrho^{}_{\ts\mu,z}\!:\U\big(t^{-1}\ggot[t^{-1}]\big)\to \U(\ggot),
\qquad X[r]\mapsto X z^r+\de_{r,-1}\ts\mu(X),\quad X\in\ggot,
\eeq
defines  a $G_\mu$-equivariant  algebra homomorphism. The image of $\z(\wh\ggot)$
under $\varrho^{}_{\ts\mu,z}$
is a commutative subalgebra $\Ac_{\mu}$ of $\U(\ggot)$ which does not depend
on $z$.
If $S\in \U\big(t^{-1}\ggot[t^{-1}]\big)$ is a homogeneous element
of degree $d$,
then regarding $\varrho^{}_{\ts\mu,z}(S)$
as a polynomial in $z^{-1}$, define the elements $S^{(i)}\in\U(\ggot)$
by the expansion
\ben
\varrho^{}_{\ts\mu,z}(S)=S^{(0)} z^{-d}+\dots+S^{(d-1)} z^{-1}+S^{(d)}.
\een
Suppose that $\mu\in\ggot^*$ is regular and that
$S_1,\dots,S_n\in \z(\wh\ggot)$ is a
complete set of homogeneous
Segal--Sugawara vectors of the respective degrees $d_1,\dots,d_n$.
Then
\begin{itemize}
\item[$\diamond$] the elements
$S^{(i)}_{k}$ with $k=1,\dots,n$ and $i=0,1,\dots,d_k-1$
are algebraically independent generators of $\Ac_{\mu}$ and
$\gr\Ac_{\mu}=\overline\Ac_{\mu}$;
\item[$\diamond$] the subalgebra $\Ac_{\mu}$ of
$\U(\ggot)$ is maximal commutative.
\end{itemize}
The first of these statements is due to \cite{fft:gm}
and the second is implied by the results of \cite{codim3}.
The elements $S^{(i)}_{k}$ generate $\Ac_{\mu}$ for any $\mu\in\ggot^*$
and the inclusion $\gr\Ac_{\mu}\supset\overline\Ac_{\mu}$ holds.
It was conjectured in \cite[Conjecture~1]{fft:gm},
that the property $\gr\Ac_{\mu}=\overline\Ac_{\mu}$ extends to
all $\mu$. Its proof in type {\sf A} was given in \cite{fm:qs}.
In what follows we give a more general argument which will imply
the conjecture in types {\sf A} and {\sf C} thus providing another proof in type {\sf A}.
We will rely on the properties of the canonical symmetrisation map
\beql{sym}
\varpi\!:\Sc(\ggot)\to\U(\ggot).
\eeq
This map was already used by Tarasov~\cite{t:cs} in type {\sf A} to show that
if $\mu\in\gl_N^*$ is semisimple, then the images of the $\mu$-shifts
of certain generators of $\Sc(\gl_N)^{\gl_N}$ under the map \eqref{sym}
generate a commutative subalgebra of $\U(\gl_N)$.
By another result of Tarasov~\cite{t:ul} this subalgebra coincides with
$\Ac_\mu$ if $\mu$ is regular semisimple.
Below we prove a similar statement for
all classical types {\sf B}, {\sf C} and {\sf D}
and all $\mu\in\ggot^*$, which allows us to
suggest that it
holds for all simple Lie algebras, see Conjecture~\ref{conj:sym} below.

In what follows we will identify any element $X\in\ggot$ with its images under the canonical
embeddings $\ggot\hra\Sc(\ggot)$ and $\ggot\hra\U(\ggot)$. It should always
be clear from the context
whether $X$ is regarded as an element of $\Sc(\ggot)$ or $\U(\ggot)$.

\subsection{Symmetrised determinants and permanents}
\label{subsec:detper}

Consider the
symmetriser $h^{(m)}$ and anti-symmetriser $a^{(m)}$ in the group algebra $\CC[\Sym_m]$
of the symmetric group $\Sym_m$. These are
the idempotents defined by
\ben
h^{(m)}=\frac{1}{m!}\ts\sum_{\si\in\Sym_m} \si
\Fand a^{(m)}=\frac{1}{m!}\ts\sum_{\si\in\Sym_m} \sgn \si\cdot \si.
\een
We let $H^{(m)}$ and $A^{(m)}$ denote their respective images under
the natural action of $\Sym_m$ on the tensor product space $(\CC^{N})^{\ot m}$.
We will denote by $P_{\si}$
the image of $\si\in\Sym_m$ with respect to this action.

For an arbitrary associative algebra $\cM$
we will identify $H^{(m)}$ and $A^{(m)}$ with the respective elements
$H^{(m)}\ot 1$ and $A^{(m)}\ot 1$ of the algebra
\beql{tenprka}
\underbrace{\End\CC^{N}\ot\dots\ot\End\CC^{N}}_m{}\ot\cM.
\eeq
Any $N\times N$ matrix $M=[M_{ij}]$ with entries in $\cM$ will be
represented as the element
\beql{mmat}
M=\sum_{i,j=1}^N e_{ij}\ot M_{ij}\in \End\CC^N\ot \cM.
\eeq
For each $a\in\{1,\dots,m\}$
introduce the element $M_a$ of the algebra
\eqref{tenprka}
by
\beql{matnota}
M_a=\sum_{i,j=1}^{N}
1^{\ot(a-1)}\ot e_{ij}\ot 1^{\ot(m-a)}\ot M_{ij}.
\eeq
The partial trace $\tr_a$
will be understood as the linear map
\ben
\tr_a\!:\End(\CC^{N})^{\ot m}\to \End(\CC^{N})^{\ot (m-1)}
\een
which acts as the trace map on the $a$-th copy of $\End\CC^N$ and
is the identity map on all the remaining copies.

For any $m=1,\dots,N$ define
the $m$-{\em th} {\em symmetrised minor} of the matrix $M$ by
\ben
\Det_m(M)=\tr_{1,\dots,m}\ts A^{(m)} M_1\dots M_m,
\een
or explicitly,
\ben
\Det_m(M)=
\frac{1}{m!}\ts \sum_{1\leqslant a_1<\dots< a_m\leqslant N}\ts
\sum_{\si,\tau\in\Sym_m}\sgn\si\tau\cdot
M_{a_{\si(1)}\tss a_{\tau(1)}}\dots M_{a_{\si(m)}\tss a_{\tau(m)}}.
\een
The {\em symmetrised determinant} of $M$ is
$\Det(M)=\Det_N(M)$. It coincides with the usual determinant
$\det(M)$ in the case of a commutative
algebra $\cM$.
For any $m\geqslant 1$
the $m$-{\em th} {\em symmetrised permanent} is defined by
\ben
\Per_m(M)=\tr_{1,\dots,m}\ts H^{(m)} M_1\dots M_m,
\een
which is written explicitly as
\ben
\Per_m(M)=\frac{1}{m!}\ts \sum_{1\leqslant a_1\leqslant\dots\leqslant a_m\leqslant N}
\frac{1}{\al_1!\dots\al_N!}\sum_{\si,\tau\in\Sym_m}
M_{a_{\si(1)}\tss a_{\tau(1)}}\dots M_{a_{\si(m)}\tss a_{\tau(m)}},
\een
where $\al_k$ denotes the multiplicity of $k\in\{1,\dots,N\}$ in the multiset
$\{a_1,\dots,a_m\}$.

Both symmetrised determinants and permanents were considered by Itoh~\cite{i:td,i:tp}
in relation to Casimir elements for classical Lie algebras.

\subsection{Generators of $\Ac_{\mu}$ in type {\sf A}}
\label{subsec:gena}

We will work with the reductive Lie algebra $\gl_N$ rather than
the simple Lie algebra $\sll_N$. We will denote the standard basis elements
of $\gl_N$ by $E_{ij}$ for $i,j=1,\dots,N$ and combine them into the $N\times N$ matrix
$E=[E_{ij}]$. Regarding the entries of $E$ as elements of the symmetric algebra $\Sc(\gl_N)$,
consider the characteristic polynomial
\ben
\det(u+E)=u^N+\Phi_1\tss u^{N-1}+\cdots+\Phi_N.
\een
Its coefficients $\Phi_1,\dots,\Phi_N$ are algebraically independent generators of
the algebra of $\gl_N$-invariants $\Sc(\gl_N)^{\gl_N}$. All coefficients $\Psi_i$ of the series
\ben
\det(1-q\tss E)^{-1}=1+\Psi_1\tss q+\Psi_2\tss q^2+\dots
\een
belong to the algebra $\Sc(\gl_N)^{\gl_N}$ and $\Psi_1,\dots,\Psi_N$
are its algebraically independent generators. Writing the matrix $E$
in the form \eqref{mmat} with $\cM=\Sc(\gl_N)$ we can represent
the generators $\Phi_m$ and $\Psi_m$ in the symmetrised form
\ben
\Phi_m=\Det_m(E) \Fand
\Psi_m=\Per_m(E).
\een
The respective $\mu$-shifts are found as the coefficients of the polynomials
in $z^{-1}$ so that
\beql{genamuant}
\Det_m(\mu+Ez^{-1})=\Phi_m\tss z^{-m}+\frac{1}{1!}\ts \di_\mu\Phi_m\tss z^{-m+1}
+\dots+\frac{1}{m!}\ts \di_\mu^m\tss\Phi_m
\eeq
and
\beql{genamuantper}
\Per_m(\mu+Ez^{-1})=\Psi_m\tss z^{-m}+\frac{1}{1!}\ts \di_\mu\Psi_m\tss z^{-m+1}
+\dots+\frac{1}{m!}\ts \di_\mu^m\tss\Psi_m
\eeq
where $\mu\in \gl_N^*$ is regarded as the
$N\times N$ matrix $\mu=[\mu_{ij}]$ with $\mu_{ij}=\mu(E_{ij})$.

The next theorem shows that
Conjecture~\ref{conj:sym} holds for each of the families $\Phi_1,\dots,\Phi_N$
and $\Psi_1,\dots,\Psi_N$.
It is the first family which was considered by Tarasov~\cite{t:cs}, \cite{t:ul}.

\bth\label{thm:conja}
Suppose that $\mu\in\gl_N^*$ is arbitrary. The algebra $\Ac_{\mu}$ is generated
by each family of
elements $\varpi(\di_\mu^{\tss k} \Phi_m)$ and $\varpi(\di_\mu^{\tss k} \Psi_m)$ with
$m=1,\dots,N$ and $k=0,1,\dots,m-1$.
\eth

\bpf
Since the coefficients
of the polynomials in \eqref{genamuant} and \eqref{genamuantper}
are already written in a symmetrised form,
their images under the symmetrisation map \eqref{sym}
are given by the same expressions $\Det_m(\mu+Ez^{-1})$ and $\Per_m(\mu+Ez^{-1})$,
where the matrix $E$ is now regarded as the element
\beql{egln}
E=\sum_{i,j=1}^N e_{ij}\ot E_{ij}\in \End\CC^N\ot \U(\gl_N).
\eeq
This follows from the
easily verified property of the symmetrisation map \eqref{sym}:
\beql{symprcy}
\varpi\!:(c_1+Y_1)\dots (c_k+Y_k)\mapsto \frac{1}{k!}\ts\sum_{\si\in\Sym_k}
(c_{\si(1)}+Y_{\si(1)})\dots (c_{\si(k)}+Y_{\si(k)})
\eeq
for any constants $c_i$ and any elements $Y_i\in\ggot$.
Therefore, we only need to show that $\Ac_{\mu}$ is generated by the coefficients
of each family of polynomials $\Det_m(\mu+Ez^{-1})$ and $\Per_m(\mu+Ez^{-1})$
with $m=1,\dots,N$, where the matrix $E$ is defined in \eqref{egln}.
However, these properties of the coefficients were already established
in \cite[Sect.~4]{fm:qs}.
\epf

\subsection{Generators of $\Ac_{\mu}$ in types {\sf B}, {\sf C} and {\sf D}}
\label{subsec:genbcd}

Define the orthogonal Lie algebras $\oa_N$ with $N=2n+1$ and $N=2n$
and symplectic Lie algebra $\spa_N$ with $N=2n$,
as subalgebras of $\gl_N$ spanned by the elements $F_{i\tss j}$ with $i,j\in\{1,\dots,N\}$,
\ben
F_{i\tss j}=E_{i\tss j}-E_{j\pr i\pr}\Fand F_{i\tss j}
=E_{i\tss j}-\ve_i\ts\ve_j\ts E_{j\pr i\pr},
\een
respectively, for $\oa_N$ and $\spa_N$.
We use the notation $i\pr=N-i+1$, and
in the symplectic case set
$\ve_i=1$ for $i=1,\dots,n$ and
$\ve_i=-1$ for $i=n+1,\dots,2n$.
We will
denote by $\ggot_N$ any of the Lie algebras $\oa_N$
or $\spa_N$.
Introduce the $N\times N$ matrix
$F=[F_{ij}]$. Regarding its entries as elements of the symmetric algebra $\Sc(\ggot_N)$,
in the symplectic case
consider the characteristic polynomial
\beql{detsp}
\det(u+F)=u^{2n}+\Phi_2\tss u^{2n-2}+\cdots+\Phi_{2n}.
\eeq
The coefficients $\Phi_2,\dots,\Phi_{2n}$ are algebraically independent generators of
the algebra of $\spa_{2n}$-invariants $\Sc(\spa_{2n})^{\spa_{2n}}$.
In the orthogonal case,
all coefficients $\Psi_{2k}$ of the series
\ben
\det(1-q\tss F)^{-1}=1+\Psi_2\tss q^2+\Psi_4\tss q^4+\dots
\een
belong to the algebra $\Sc(\oa_N)^{\oa_N}$.
In the case of even $N=2n$
we also define the {\em Pfaffian\/}
by
\ben
\Pf\tss F
=\frac{1}{2^nn!}\sum_{\si\in\Sym_{2n}}\sgn\si\cdot
F_{\si(1)\ts\si(2)'}\dots
F_{\si(2n-1)\ts\si(2n)'}.
\een
The coefficients $\Psi_2,\dots,\Psi_{2n}$
are algebraically independent generators of the algebra
$\Sc(\oa_{2n+1})^{\oa_{2n+1}}$, while
the elements $\Psi_2,\dots,\Psi_{2n-2},\Pf\tss F$
are algebraically independent generators of $\Sc(\oa_{2n})^{\oa_{2n}}$.

Write the matrix $F$
in the form \eqref{mmat} with $\cM=\Sc(\ggot_N)$ and represent
the generators $\Phi_m$ and $\Psi_m$ for even values of $m$ in the symmetrised form
\ben
\Phi_m=\Det_m(F) \Fand
\Psi_m=\Per_m(F).
\een
Then the respective $\mu$-shifts are found as the coefficients of the polynomials
in $z^{-1}$,
\beql{genamuantf}
\Det_m(\mu+Fz^{-1})=\Phi_m\tss z^{-m}+\frac{1}{1!}\ts \di_\mu\Phi_m\tss z^{-m+1}
+\dots+\frac{1}{m!}\ts \di_\mu^m\tss\Phi_m
\eeq
and
\beql{genamuantperf}
\Per_m(\mu+Fz^{-1})=\Psi_m\tss z^{-m}+\frac{1}{1!}\ts \di_\mu\Psi_m\tss z^{-m+1}
+\dots+\frac{1}{m!}\ts \di_\mu^m\tss\Psi_m,
\eeq
where $\mu\in \ggot_N^*$ is regarded as the
$N\times N$ matrix $\mu=[\mu_{ij}]$ with $\mu_{ij}=\mu(F_{ij})$.
The $\mu$-shifts of the Pfaffian $\Pf\tss F$
are the coefficients of the polynomial
\ben
\Pf\ts \big(\mu+F\tss z^{-1}\big)=\pi^{}_{(0)}z^{-n}+
\dots+\pi^{}_{(n-1)}z^{-1}
+\pi^{}_{(n)},\qquad \pi^{}_{(k)}\in\Sc(\oa_{2n}),
\een
where
\beql{pfaffea}
\Pf\ts \big(\mu+F\tss z^{-1}\big)
=\frac{1}{2^nn!}\sum_{\si\in\Sym_{2n}}\sgn\si\cdot
\big(\mu+F\tss z^{-1}\big)_{\si(1)\ts\si(2)'}\dots
\big(\mu+F\tss z^{-1}\big)_{\si(2n-1)\ts\si(2n)'}.
\eeq

The following theorem implies
Conjecture~\ref{conj:sym} for the orthogonal and symplectic Lie algebras
for the families of generators of $\Sc(\ggot_N)^{\ggot_N}$
described above.

\bth\label{thm:conjbcd}
Suppose that $\mu\in\ggot_N^*$ is arbitrary.
\begin{itemize}
\item[$\diamond$] The family
$\varpi(\di_\mu^{\tss p} \Phi_m)$ with
$m=2,4,\dots,2n$ and $p=0,1,\dots,m-1$
generates the algebra $\Ac_{\mu}$
in type {\sf C}.
\item[$\diamond$] The family
$\varpi(\di_\mu^{\tss p} \Psi_m)$ with
$m=2,4,\dots,2n$ and $p=0,1,\dots,m-1$
generates the algebra $\Ac_{\mu}$
in type {\sf B}.
\item[$\diamond$] The family
$\varpi(\di_\mu^{\tss p} \Psi_m)$ with
$m=2,4,\dots,2n-2$ and $p=0,1,\dots,m-1$
together with the elements $\varpi(\pi^{}_{(k)})$ for $k=0,\dots,n-1$
generate the algebra $\Ac_{\mu}$
in type {\sf D}.
\end{itemize}
\eth

\bpf
It follows from the results of \cite{m:ff} that the generators
of $\Sc(\ggot_N)^{\ggot_N}$ introduced above admit the well-defined forms
\ben
\Psi_{m}=\ga_{m}(N)\ts\tr^{}_{1,\dots,m}\ts S^{(m)}
F_1\dots F_{m} \Fand
\Phi_{m}=\ga_{m}(-2n)\ts\tr^{}_{1,\dots,m}\ts S^{(m)}
F_1\dots F_{m}
\een
in the orthogonal and symplectic case, respectively, where
\ben
\ga_m(\om)=\frac{\om+m-2}{\om+2\tss m-2},
\een
and $S^{(m)}\in \End(\CC^N)^{\ot\ts m}$ is the {\em Brauer algebra symmetriser}
which admits a few equivalent expressions; see also \cite{m:so}.
The corresponding $\mu$-shifts are found as the coefficients of the polynomial
\beql{genamuantbcd}
\ga_{m}(\om)\ts\tr^{}_{1,\dots,m}\ts S^{(m)}
\big(\mu^{}_1+F_1\tss z^{-1}\big)\dots
\big(\mu^{}_{m}+F_{m}\tss z^{-1}\big),
\eeq
where we extend
notation \eqref{matnota} to the matrix $\mu$
and assume the specialisations
$\om=N$ and $\om=-2n$ in the orthogonal and symplectic case, respectively.
Due to \eqref{symprcy},
the images of the polynomials \eqref{genamuantf} and \eqref{genamuantperf}
under the symmetrisation map \eqref{sym}
are given by the same expressions $\Det_m(\mu+Fz^{-1})$ and $\Per_m(\mu+Fz^{-1})$
where the matrix $F$ is now regarded as the element
\beql{fbcd}
F=\sum_{i,j=1}^N e_{ij}\ot E_{ij}\in \End\CC^N\ot \U(\ggot_N).
\eeq
The same observation shows that
the image of the polynomial \eqref{genamuantbcd}
under the symmetrisation map \eqref{sym}
is given by the expression \eqref{genamuantbcd}
where the matrix $F$ is given by
\eqref{fbcd}. This allows us to conclude that
the following identities hold for polynomials with coefficients in $\U(\ggot_N)$:
\beql{detsym}
\Det_m(\mu+Fz^{-1})=\ga_{m}(-2n)\ts\tr^{}_{1,\dots,m}\ts S^{(m)}
\big(\mu^{}_1+F_1\tss z^{-1}\big)\dots
\big(\mu^{}_{m}+F_{m}\tss z^{-1}\big)
\eeq
in the symplectic case, and
\beql{persym}
\Per_m(\mu+Fz^{-1})=\ga_{m}(N)\ts\tr^{}_{1,\dots,m}\ts S^{(m)}
\big(\mu^{}_1+F_1\tss z^{-1}\big)\dots
\big(\mu^{}_{m}+F_{m}\tss z^{-1}\big)
\eeq
in the orthogonal case.

Now consider the orthogonal and symplectic Lie algebras simultaneously and
define the polynomials $\phi^{}_{m}(z)$ in $z^{-1}$
by
\beql{expagasm}
\phi^{}_{m}(z)=\ga_m(\om)\ts\tr^{}_{1,\dots,m}\ts S^{(m)}
\big({-}\di_z+\mu^{}_1+F_1\tss z^{-1}\big)\dots
\big({-}\di_z+\mu^{}_m+F_m\tss z^{-1}\big)\ts 1,
\eeq
where $\di_z$ is understood
as the differential operator with $\di_z\ts 1=0$
so that
\ben
\phi^{}_{m}(z)=\phi^{}_{m\ts(0)}z^{-m}
+\dots+\phi^{}_{m\ts(m-1)}z^{-1}
+\phi^{}_{m\ts(m)},\qquad \phi^{}_{m\ts(k)}\in \U(\ggot_N).
\een
As for the expression \eqref{genamuantbcd},
the right hand side of \eqref{expagasm} in the symplectic case is assumed to be written
in a certain equivalent form which is well-defined for all $1\leqslant m\leqslant 2n+1$;
see \cite{m:ff} and \cite{m:so} for details. The same assumption will apply to all
expressions of this kind throughout the rest of the proof.

In the case $\ggot_N=\oa_{2n}$ the $\varpi$-image of $\Pf\ts \big(\mu+F\tss z^{-1}\big)$
coincides with the expression defined by \eqref{pfaffea}, where $F$ is now defined by
\eqref{fbcd}.
With this interpretation of $F$ we can write
\beql{pfaffnonc}
\Pf\ts \big(\mu+F\tss z^{-1}\big)=\varpi(\pi^{}_{(0)})z^{-n}+
\dots+\varpi(\pi^{}_{(n-1)})z^{-1}
+\varpi(\pi^{}_{(n)}).
\eeq

By the general properties of the algebra $\Ac_{\mu}$ from \cite{fft:gm} which we
recalled in the beginning of this
section and the results of \cite{m:ff},
given any $\mu\in\ggot^*_N$, all coefficients of the polynomials
$\phi^{}_{m}(z)$
belong to the commutative subalgebra $\Ac_{\mu}$ of $\U(\ggot_N)$.
All coefficients of the polynomial \eqref{pfaffnonc}
belong to the subalgebra $\Ac_{\mu}$ of $\U(\oa_{2n})$.
Moreover,
the elements
\ben
\phi^{}_{2k\ts(p)}\qquad
\text{with}\quad k=1,\dots,n \fand p=0,1,\dots,2k-1
\een
generate the algebra $\Ac_{\mu}$ in the cases {\sf B} and {\sf C},
while the elements
\ben
\phi^{}_{2k\ts(p)}\qquad
\text{with}\quad k=1,\dots,n-1 \fand p=0,1,\dots,2k-1
\een
together with
$\varpi(\pi^{}_{(0)}),\dots,\varpi(\pi^{}_{(n-1)})$
generate the algebra $\Ac_{\mu}$ in the case
{\sf D}.

Due to \eqref{expagasm},
$\phi^{}_{m\tss (p)}$ is found as the coefficient of $z^{-m+p}$ in the expression
\ben
\sum_{i_1<\dots<i_p}
\sum_{j_1<\dots<j_{m-p}}
\ga_m(\om)\ts\tr_{1,\dots,m}\ts S^{(m)} \mu^{}_{i_1}\dots \mu^{}_{i_p}
\big({-}\di_z+F^{}_{j_1}z^{-1}\big)\dots
\big({-}\di_z+F^{}_{j_{m-p}}z^{-1}\big)\ts 1,
\een
summed over disjoint subsets of indices $\{i_1,\dots,i_p\}$ and
$\{j_1,\dots,j_{m-p}\}$ of $\{1,\dots,m\}$.
Note that for any $\si\in\Sym_m$ we have $S^{(m)}P_{\si}=P_{\si}\tss S^{(m)}=\pm\tss S^{(m)}$.
Hence,
applying the cyclic property of trace together with the relations $P_{\si}\tss\mu_i=\mu_{\si(i)}P_{\si}$
and $P_{\si}F_j=F_{\si(j)}P_{\si}$, we get
\ben
\phi^{}_{m\tss (p)}=z^{m-p}\ts\binom{m}{p}\ts\ga_m(\om)\ts\tr_{1,\dots,m}
\ts S^{(m)} \mu^{}_1\dots \mu^{}_p
\big({-}\di_z+F^{}_{p+1}z^{-1}\big)\dots
\big({-}\di_z+F^{}_mz^{-1}\big)\ts 1.
\een
Furthermore,
\ben
\bal
z^{m-p}\ts\big({-}\di_z+F^{}_{p+1}z^{-1}\big)&\dots
\big({-}\di_z+F^{}_mz^{-1}\big)\ts 1\\[0.3em]
{}&=F^{}_{p+1}\dots F^{}_m+\text{a linear combination of }\  \{F^{}_{a_1}\dots F^{}_{a_s}\}
\eal
\een
with $p+1\le a_1<\dots<a_s\le m$. Applying again the cyclic property of trace and
appropriate conjugations by the elements $P_{\si}$, we bring the expression for
$\phi^{}_{m\tss (p)}$ to the form
\begin{multline}
\non
\phi^{}_{m\tss (p)}=\binom{m}{p}\ts\ga_m(\om)\ts\tr_{1,\dots,m}\ts S^{(m)} \mu^{}_1\dots \mu^{}_p
\ts F^{}_{p+1}\dots F^{}_m\\
{}+\sum_{r=p+1}^{m-1} c_r\ts\ga_m(\om)\ts\tr_{1,\dots,m}\ts
S^{(m)} \mu^{}_1\dots \mu^{}_p
\ts F^{}_{p+1}\dots F^{}_r
\end{multline}
for certain constants $c_r$.
For the partial trace of the symmetriser we have
\ben
\tr_m\ts\ga_m(\om)\ts S^{(m)}=\pm\tss\frac{\om+m-2}{m}\ts
\ga_{m-1}(\om)\ts S^{(m-1)},
\een
with the plus and minus sign taken in the
orthogonal and symplectic case, respectively (assuming $m\leqslant n$ for the latter);
see \cite{m:ff}. This yields
\begin{multline}
\non
\phi^{}_{m\tss (p)}=\binom{m}{p}\ts\ga_m(\om)\ts\tr_{1,\dots,m}\ts S^{(m)} \mu^{}_1\dots \mu^{}_p
\ts F^{}_{p+1}\dots F^{}_m\\
{}+\sum_{r=p+1}^{m-1} d_r\ts\ga_r(\om)\ts\tr_{1,\dots,r}\ts
S^{(r)} \mu^{}_1\dots \mu^{}_p
\ts F^{}_{p+1}\dots F^{}_r
\end{multline}
for certain constants $d_r$.
Introduce the coefficients of the respective polynomials
on the right hand sides of
\eqref{detsym} and \eqref{persym} by
\ben
\ga_{m}(\om)\ts\tr^{}_{1,\dots,m}\ts S^{(m)}
\big(\mu^{}_1+F_1\tss z^{-1}\big)\dots
\big(\mu^{}_{m}+F_{m}\tss z^{-1}\big)
=\vp^{}_{m\ts(0)}z^{-m}
+\dots+\vp^{}_{m\ts(m-1)}z^{-1}
+\vp^{}_{m\ts(m)}
\een
with $\vp^{}_{m\tss(p)}\in \U(\ggot_N)$.
The same argument as for the coefficients $\phi^{}_{m\tss (p)}$
gives
\ben
\vp^{}_{m\tss (p)}=\binom{m}{p}\ts\ga_m(\om)\ts\tr_{1,\dots,m}\ts S^{(m)} \mu^{}_1\dots \mu^{}_p
\ts F^{}_{p+1}\dots F^{}_m.
\een
This yields a triangular system of linear
relations
\ben
\phi^{}_{m\tss (p)}=\vp^{}_{m\tss (p)}+\sum_{r=p+1}^{m-1} d\pr_r\ts \vp^{}_{r\ts (p)}.
\een
Therefore, we may now conclude that the elements $\vp^{}_{m\tss (p)}$
with even $m=2,4,\dots,2n$ and $p=0,1,\dots,m-1$
generate the algebra $\Ac_{\mu}$ in types {\sf B} and {\sf C}, while the elements $\vp^{}_{m\tss (p)}$
with $m=2,4,\dots,2n-2$ and $p=0,1,\dots,m-1$ together with
$\varpi(\pi^{}_{(0)}),\dots,\varpi(\pi^{}_{(n-1)})$
generate the algebra $\Ac_{\mu}$ in type
{\sf D}. The proof is completed by taking into account the relations \eqref{detsym}
and \eqref{persym}.
\epf

\bcj\label{conj:sym}
For any simple Lie algebra $\ggot$, 
there exist generators $H_1,\dots,H_n$ of the algebra
$\Sc(\ggot)^{\ggot}$ such that
for any element $\mu\in\ggot^*$, 
the images $\varpi(\partial_\mu^k H_i)$ of their $\mu$-shifts
with respect to the symmetrisation map $\varpi$ generate the algebra $\Ac_{\mu}$.
\ecj

\section{Quantisations of MF subalgebras with a non-regular $\mu$}

Take any $\mu\in\gt g^*$.
As we already know, $\gr\Ac_\mu$ is a Poisson-commutative subalgebra of $\cS(\gt g)$ and
$\oA_\mu\subset\gr\Ac_\mu$ by \cite[Prop.~3.12]{fft:gm}.
As can be seen from the construction \eqref{evalr}, $\Ac_\mu\subset{\mathcal U}(\gt g)^{G_\mu}$
and thereby $\gr\Ac_\mu\subset\cS(\gt g)^{G_\mu}$.
According to Proposition~\ref{l-tr.deg},
${\rm tr.deg}(\gr\Ac_\mu)\le \frac{1}{2}\dim (G\mu)+\rk\gt g$. At the same time,
${\rm tr.deg}\,\oA_\mu = \frac{1}{2}\dim (G\mu)+\rk\gt g$ by Lemma~\ref{lm-bol}.  Hence we have the following general result.

\begin{prop}
For any reductive $\gt g$ and any $\mu\in\gt g^*$, $\gr\Ac_\mu$ is an algebraic extension of $\oA_\mu$.
\qed \end{prop}

Suppose that $\oA_\mu$ is a maximal Poisson-commutative subalgebra of $\cS(\gt g)^{\gt g_\mu}$.
Then necessary $\gr{\mathcal A}_\mu =\oA_\mu$.  In view of Corollary~\ref{A-C},
the FFTL-conjecture in types {\sf A} and {\sf C} follows.

\begin{thm}\label{thm-AC}
Suppose that $\gt g$ is of type {\sf A} or {\sf C}. Then  $\gr{\mathcal A}_\mu =\oA_\mu$ for each
$\mu\in\gt g^*$. \qed
\end{thm}

We can  rely
on Proposition~\ref{prop-cone} to
conclude that $\gr{\mathcal A}_\mu =\oA_\mu$ for some $\mu\in\gt g^*\cong\gt g$ in the other types.

\begin{ex}[The minimal nilpotent orbit]
Let $\gamma\in\gt g$ be a minimal nilpotent element in a simple Lie algebra
$\gt g$. Suppose that
$\gt g$ is not of type {\sf E}$_8$. Then there is a g.g.s.
for $\gamma$ and $\gt g_\gamma$ has the codim-2 property, see \cite{ppy}.
Hence $\gr{\mathcal A}_\gamma =\oA_\gamma$.
\end{ex}

\begin{ex}[The subregular case]
Keep the assumption that $\gt g$ is simple and assume that $\dim(G\mu)=\dim\gt g-\rk\gt g-2$.
Let $\gamma$ be as in Lemma~\ref{curve}. Then $G\gamma$ is the subregular nilpotent orbit.
There is a g.g.s. for $\gamma$ and $\gt g_\gamma$ has the codim-2 property if $\gt g$ is not of type
${\sf G}_2$, see \cite[Sect.~6]{par-contr}.  In view of Proposition~\ref{prop-cone}, we have
$\gr{\mathcal A}_\mu =\oA_\mu$ outside of type {\sf G}$_2$.
\end{ex}

In types {\sf A} and {\sf C}, it is possible to describe the generators of ${\mathcal A}_\mu$ explicitly.
If $\gt g=\gt{gl}_N$, then $\Phi_1, \ldots,\Phi_N$ is a g.g.s. for any nilpotent $\gamma\in\gt g$;
if $\gt g=\gt{sp}_{2n}$, then $\Phi_2,\ldots,\Phi_{2n}$ is a g.g.s. for any nilpotent $\gamma\in\gt g$,
see \cite{ppy}. The degrees of ${^\gamma\Phi_i}$ can be found in \cite[Sect.~4]{ppy}.
We give  more details below.

\subsubsection*{Type {\sf A}}
Suppose that $\gt  g=\gl_N$ and that $\mu\in\gl_N$ has the Jordan blocks as in
\eqref{jord-s}. Then
the nilpotent element of Lemma~\ref{curve} is given by the
partition $\Pi$ defined in \eqref{jord-n}. Let
$\Pi=(\beta_1,\ldots,\beta_M)$.  Using the results
of Sections~\ref{sec-MF} and \ref{subsec:gena} we find that the algebra
${\mathcal A}_\mu$ is freely generated
by the elements
$\varpi(\partial^k_\mu\Phi_m)$ with $1\le m\le N$ and $0\le k\le m{-}{\bf r}(m)$, where
${\bf r}(m)$ is uniquely determined by the conditions
\ben
\sum_{i=1}^{{\bf r}(m){-}1} \beta_i < m \le \sum_{i=1}^{{\bf r}(m)} \beta_i.
\een
To give an equivalent definition of ${\bf r}(m)$,
consider the row-tableau of shape $\Pi$ which is obtained by writing
the numbers $1,2,\dots,N$ consecutively from left to right in the boxes of each row
of the Young diagram $\Pi$
beginning from the top row. Then ${\bf r}(m)$ equals the row number of $m$ in the tableau.

Note that
in the case where $\Pi$ corresponds to a nilpotent element $\mu$,
the elements  $\varpi(\partial^k_\mu\Phi_m)$ with $k>m{-}{\bf r}(m)$ are equal to zero.
In the general case,
associate the elements of the family $\Phi_{m\tss k}=\varpi(\partial^k_\mu\Phi_m)$
with the boxes of the diagram $\Ga=(N,N-1,\dots,1)$,
as illustrated:
\ben
\Ga\quad=\quad
\begin{matrix}
\Phi_{N\tss N-1} & \Phi_{N\tss N-2} & \dots & \Phi_{N\tss 1} & \Phi_{N\ts 0}\\[0.5em]
\quad \Phi_{N-1\tss N-2} & \quad\Phi_{N-1\tss N-3} & \dots & \quad\Phi_{N-1\ts 0} & \\
\dots & \dots & \dots & & \\[0.5em]
\Phi_{\tss 2\ts 1}\phantom{-,} & \Phi_{\tss 2\ts 0}\phantom{-,} & & &\\[0.5em]
\Phi_{\tss 1\ts 0}\phantom{-,} & & & &
\end{matrix}
\een
Then the free generators of $\Ac_{\mu}$ correspond to the skew diagram
$\Ga/\si$, where
\ben
\si=\big({\bf r}(N)-1,\dots,{\bf r}(1)-1\big).
\een
This agrees with the results of \cite{fm:qs} which were applied to a few different
families of generators, where the transpose of $\Ga/\si$
was used instead. The transpose of $\si$ is found by
\ben
\si\pr=\big(\beta_2+\dots+\beta_M,\beta_3+\dots+\beta_M,\dots,\beta_M\big)
\een
which
coincides with the diagram $\ga$
in the notation of that paper. It was shown there that any generator
$\Phi_{m\tss k}$ associated with a box of the diagram $\si$ is represented
by a {\em linear combination} of the generators corresponding to the boxes
of $\Ga/\si$ which occur in the column containing $\Phi_{m\tss k}$.

\subsubsection*{Type {\sf C}}
Now suppose that $\gt g=\spa_{2n}$. Keep the notation of \eqref{jord-s} for the Jordan blocks of
$\mu\in\spa_{2n}$. The nilpotent element of Lemma~\ref{curve} is associated with the partition
$\Pi_{\gamma}=(\beta_1,\ldots,\beta_M)$
produced in Section~\ref{subsec:shlim}.
Similar to type {\sf A}, for each $m=1,\dots,n$
define positive integers ${\bf r}(2m)$ by the conditions
\ben
\sum_{i=1}^{{\bf r}(2m){-}1} \beta_i < 2m \le \sum_{i=1}^{{\bf r}(2m)} \beta_i.
\een
Equivalently, they can be defined in terms of the initial Young diagram $\Pi$ given
in \eqref{jord-n} by the following rule.
Consider the row-tableau of shape $\Pi$ which is obtained by writing
the numbers $1,2,\dots,2n$ consecutively from left to right in the boxes of each row
of $\Pi$
beginning from the top row. Then ${\bf r}(2m)$ equals the
row number of $2m$ in the tableau.

The algebra ${\mathcal A}_\mu$ is freely generated by
the elements $\Phi_{2m\ts k}=\varpi(\partial_\mu^k\Phi_{2m})$
with $1\le m\le n$ and $0\le k\le 2m-{\bf r}(2m)$.
To illustrate this result in terms of diagrams,
associate the elements of the family $\Phi_{2m\ts k}$
with the boxes of the diagram $\Ga=(2n,2n-2,\dots,2)$,
as follows:
\ben
\Ga\quad=\quad
\begin{matrix}
\Phi_{2n\tss 2n-1} & \Phi_{2n\tss 2n-2} & \dots & \Phi_{2n\tss 2} & \Phi_{2n\ts 1}& \Phi_{2n\ts 0}\\[0.5em]
\quad \Phi_{2n-2\ts 2n-3} & \quad\Phi_{2n-2\ts 2n-4} & \dots & \quad\Phi_{2n-2\ts 0} & &\\
\dots & \dots & \dots & & &\\[0.5em]
\Phi_{\tss 2\ts 1}\phantom{-,} & \Phi_{\tss 2\ts 0}\phantom{-,} & & & &
\end{matrix}
\een

\bigskip
\noindent
Then the free generators of $\Ac_{\mu}$ correspond to the skew diagram
$\Ga/\si$, where
\ben
\si=\big({\bf r}(2n)-1,{\bf r}(2n-2)-1,\dots,{\bf r}(2)-1\big).
\een
If $\mu$ is regular then $\si$ is empty and all generators in $\Ga$
are algebraically independent. In the other extreme case where $\mu=0$
the diagram $\si$ is $(2n-1,2n-3,\dots,1)$ and $\Ac_{\mu}$ coincides
with the center of the universal enveloping algebra $\U(\spa_{2n})$.
It is freely generated by the elements $\Phi_{2\ts 0},\Phi_{4\ts 0},\dots,\Phi_{2n\ts 0}$.

For another illustration consider $\spa_{10}$ and suppose that
$\mu\in \spa_{10}$ has the zero eigenvalue with the corresponding Young
diagram $(1,1)$ and two opposite sign eigenvalues, each corresponding
to the diagram $(2,1,1)$. Then $\Pi=(5,3,2)$ with
\ben
{\bf r}(2)={\bf r}(4)=1,\qquad {\bf r}(6)={\bf r}(8)=2\Fand {\bf r}(10)=3.
\een
Hence $\si=(2,1,1,0,0)$ and the algebra $\Ac_{\mu}$ is freely generated by
all elements in $\Ga$ except for $\Phi_{10\ 9},\Phi_{10\ 8},\Phi_{8\ts 7}$
and $\Phi_{6\ts 5}$.

\begin{rmk}
Let $\gamma\in\gt g^*$ be a nilpotent element and assume that there is a g.g.s.
$H_1,\ldots,H_n$ for $\gamma$. Then $\oA_\gamma$ is freely generated
by $F_1,\ldots,F_s$\,, where each $F_j$ is a $\gamma$-shift $\partial_\gamma^k H_i$\,, see Section~\ref{sec-MF}. Take  elements ${\mathcal F}_j\in\Ac_\gamma$ such that the symbol of
${\mathcal F}_j$ is $F_j$. Then the commutative subalgebra
$\widetilde{\Ac}_\gamma\subset\U(\gt g)$
generated by $\Fc_1,\dots,\Fc_s$
is a
quantisation of $\oA_\gamma$, i.e.,
$\gr \widetilde{\Ac}_\gamma=\oA_\gamma$,  and therefore solves Vinberg's problem.
However we cannot claim that $\widetilde{\Ac}_\gamma = \Ac_\gamma$.
\end{rmk}

Suppose now that $\gt g=\gt{o}_N$. There are nilpotent elements $\gamma\in\gt g$ that have a g.g.s. and
the codim-2 property \cite[Thm~4.7]{ppy}. There are some other
elements that have only a g.g.s.
\cite[Lemmas~4.5,\,4.6]{ppy}.  We postpone the detailed exploration of subalgebras $\Ac_\mu\subset\U(\gt g)$ until a forthcoming paper.


\section{Quantisations of Vinberg's limits}
Take $\gamma\in\gt g^*$,  $\mu\in\gt g^*_{\rm reg}$
and let $\oC_{\gamma,\mu}$ be Vinberg's limit  at $\gamma$ along $\mu$
as defined in Section~\ref{sec-limits}. By the construction \eqref{evalr}, we have
${\mathcal A}_{\gamma{+}u\mu}\subset{\mathcal U}(\gt g)[u]$.
Let
$$
{\mathcal C}_{\gamma,\mu}=\lim\limits_{u\to 0} {\mathcal A}_{\gamma{+}u\mu}
$$
be the limit taken in the same sense as in \eqref{V-limit-def}.  
Clearly ${\mathcal C}_{\gamma,\mu}$ is a commutative subalgebra.
One could expect that this subalgebra is a quantisation of $\oC_{\gamma,\mu}$
subject to some reasonable conditions. However, this is not necessarily the case because
the operations of taking the limit and the symbol may not commute. 


\begin{lm}\label{lim-sim}
Suppose that $\oA_{\gamma+u\mu}$ is freely generated by
some elements
$F_1(u),\ldots,F_s(u)$, depending on $u$, such that the nonzero vectors
  $F_j\in \lim\limits_{u\to 0} \mK F_j(u)$ are algebraically independent. Suppose further  that
${\mathcal A}_{\gamma+u\mu}$ is generated by ${\mathcal F}_j(u)=\varpi\big((F_j(u)\big)$.
Then $\gr{\mathcal C}_{\gamma,\mu}=\oC_{\gamma,\mu}$.
\end{lm}
\begin{proof}
Since ${\mathcal F}_j(u)$ is the symmetrisation of $F_j(u)$, we have
$$\lim\limits_{u\to 0} \mK {\mathcal F_j}(u)=
\varpi\big(\lim\limits_{u\to 0} \mK F_j(u)\big)=\varpi(\mK F_j)=\mK\varpi(F_j).$$
Furthermore, the elements $\varpi(F_j)$ are algebraically independent.
Hence ${\mathcal C}_{\gamma,\mu}$ is freely generated by $\varpi(F_j)$.
Because $\oC_{\gamma,\mu}$ is (freely) generated by $F_j$,    the result follows.
\end{proof}

\begin{prop}\label{cor:limvin}
Let $\gt g$ be of type {\sf A} or {\sf C}. Suppose that
the element $\gamma\in\gt g^*\cong\gt g$ is nilpotent,
$\mu\in\gt g^*_{\rm reg}$ and $\bar\mu:=\mu|_{\gt g_\gamma}\in(\gt g_\gamma^*)_{\rm reg}$.
Then $\gr{\mathcal C}_{\gamma,\mu}=\oC_{\gamma,\mu}$.
\end{prop}
\begin{proof} Under the assumptions,  $\gt g_\gamma$ has
the codim-2 property \cite{ppy}. Furthermore,
let $\{H_1,\ldots,H_n\}\subset\cS(\gt g)^{\gt g}$ be the set of generators, where
$H_i=\Phi_{i}$ in type {\sf A} and $H_i=\Phi_{2i}$ in type {\sf C}.
Then $H_1,\ldots,H_n$ is a g.g.s. for any nilpotent $\gamma\in\gt g$ \cite{ppy}.
Since $\bar\mu\in(\gt g_\gamma)^*$, Theorem~\ref{VL-nilp} applies and asserts that
$\oC_{\gamma,\mu}$ has a set of algebraically independent generators, say
$F_1,\ldots,F_{{\bf b}(\gt g)}$. Here
$F_j\in\lim\limits_{u\to 0} \mK F_j(u)$, where $F_j(u)=\partial_{\gamma{+}u\mu}^k H_i$.
By Theorems~\ref{thm:conja} and \ref{thm:conjbcd}, each subalgebra ${\mathcal A}_{\gamma{+}u\mu}$
is generated by $\varpi\big(F_j(u)\big)$.  Hence Lemma~\ref{lim-sim} applies thus completing the proof.
\end{proof}



\subsection{Limits along regular series and symmetrisation} Limits of Mishchenko--Fomenko subalgebras have been studied since \cite{v:sc}.  A rather general definition and a detailed discussion can be found in
\cite{V14}. The following construction will be sufficient for our purposes.

Let $\gt h\subset\gt g$ be a Cartan subalgebra. Consider a sequence of elements $h(0),\ldots,h(\ell)\in\gt h$ such that
$\gt g_{h(0)}\cap\gt g_{h(1)}\cap\ldots\cap\gt g_{h(\ell)}=\gt h$ for the centralisers of $h(m)$.
Set $$\mu(u)=h(0)+uh(1)+u^2h(2)+\dots + u^\ell h(\ell).$$ Then
$\oA_{\mu(u)}\subset\cS(\gt g)[u]$. Further, set
$$
\oC=\lim\limits_{u\to 0}\oA_{\mu(u)}
$$
in the same sense as in \eqref{V-limit-def}.
In our previous considerations $\ell$ was equal to $1$, $h(1)$ was regular, but neither of
$h(0)$, $h(1)$ had to be semisimple.

The first interesting property is that ${\rm tr.deg}\,\oC={\bf b}(\gt g)$ \cite{v:sc,bk-GK,V14}. Another one is that $\oC$ is a free algebra \cite{Vitya}.
Moreover, each $\oC$ is a maximal Poisson-commutative subalgebra of
$\cS(\gt g)$ \cite{t:max}. 
In type {\sf A},  the symmetrisation map $\varpi$ provides a quantisation of
$\oC$ and a   quantisation of this subalgebra  is 
unique 
\cite{t:cs,t:ul}.

Set $\gt g_0=\gt g_{h(0)}$ and $\gt g_i=\gt g_{i{-}1}\cap \gt g_{h(i)}$
for each $i\ge 1$.
Regard $h(i)$ as a linear function on $\gt g$ and on each
$\gt g_j$. Denote by $\oA^{\tss(i)}_{h(i)}$ the Mishchenko--Fomenko subalgebra
of $\cS(\gt g_{i{-}1})$ associated with $h(i)$.

\begin{thm}[{\cite[Thm~1]{Vitya}}]
The  algebra $\oC$ is generated by the Mishchenko--Fomenko subalgebras
$\oA_{h(0)}$ and
$\oA^{\tss(i)}_{h(i)}$ with $1\le i\le \ell$, and by $\gt h=\gt g_\ell$.
\end{thm}

Recall that the subalgebra $\Ac_{h(0)}$ of $\U(\ggot)$ is defined
as the image of the Feigin--Frenkel centre under the
homomorphism \eqref{evalr} with $\mu=h(0)$. Consider the subalgebras
$\Ac^{(i)}_{h(i)}\subset\Uc(\gt g_{i{-}1})$ for $i\ge 1$
defined in the same way.
Let $\widetilde{\mathcal C}\subset {\mathcal U}(\gt g)$ be the subalgebra generated by
$\Ac_{h(0)}$ and $\Ac^{(i)}_{h(i)}\subset\Uc(\gt g_{i{-}1})$ with $1\le i\le \ell$, and by $\gt h$.
Since $\Ac^{(i)}_{h(i)}$ commutes with $\gt g_i$, the subalgebra
$\widetilde{\mathcal C}$ is commutative.

\begin{thm}
Let $\gt g$ be a reductive Lie algebra. Then $\gr\widetilde{\mathcal C}=\oC$.
\end{thm}
\begin{proof}
By the construction $\oC\subset \gr\widetilde{\mathcal C}$. Since $\oC$ is a maximal Poisson-commutative
subalgebra \cite{t:max}, we have the equality here.
\end{proof}

According to \cite[Sect.~10]{cris}, we have $\widetilde{\mathcal C}={\mathcal C}$, where 
${\mathcal C}=\lim\limits_{u\to 0}\Ac_{\mu(u)}\subset\U(\gt g)$. 

Suppose now that $\gt g$ is a classical Lie algebra. Then each $\gt g_i$ is also
a classical Lie algebra and every  simple factor of $\gt g_i$ is either of type ${\sf A}$ or of the same type as $\gt g$.
In this case, $\oC$ has a set of generators $F_1,\ldots,F_{\bf b(\gt g)}$ such that
$\widetilde{\mathcal C}$ is generated by $\varpi(F_i)$, see Section~\ref{sec:qs}.

\begin{ex}[{cf. \cite[Lemma~4]{r:si}}] Take $\mu(u)=E_{11}+uE_{22}+\dots+u^{N{-}1}E_{NN}\in\gt{gl}_N$.
Then $\oC=\lim\limits_{u\to 0} \oA_{\mu(u)}$ is generated by $\gt h$, $\cS(\gt g)^{\gt g}$, and
$\cS(\gt g_i)^{\gt g_i}$ with $0\le i\le N{-}2$ \cite{v:sc}. It is the graded image  algebra of
the Gelfand--Tsetlin subalgebra ${\mathcal{GT}}(\gl_N)\subset\U(\gt{gl}_N)$.
For $\gamma=E_{11}$,
we have $\deg{^\gamma\Phi_m}=m{-}1$ and $\di_{\ga}\Phi_1=1$ and so
\ben
\oA_{\ga}=\CC[\Phi_1,\dots,\Phi_N,\di_{\ga}\Phi_2,\dots,\di_{\ga}\Phi_N].
\een
Arguing as in the proof of Theorem~\ref{VL-nilp}, we conclude that
$\oC$ is freely generated by 
the lowest $u$-components of the elements
$\partial_{\mu(u)}^k \Phi_m$ with $1\le m\le N$ and $0\le k\le m{-}1$.
At the same time, $\Ac_{\mu(u)}$ is generated by $\varpi(\partial_{\mu(u)}^k \Phi_m)$.
Therefore, applying Lemma~\ref{lim-sim} we get $\gr{\mathcal C}=\oC$. It is clear
that ${\mathcal C}={\mathcal{GT}}(\gl_N)$.
\qed
\end{ex}

The main property of ${\mathcal{GT}}(\gl_N)$ is that this subalgebra has a simple spectrum in any
irreducible finite-dimensional $\GL_N$-module $V$. It was noticed in \cite{r:si} that
$\Ac_\mu\subset\U(\gt{gl}_N)$ has the same property if $\mu\in\gt h$ is sufficiently generic.
The observation was extended to all reductive Lie algebras $\gt g$ in \cite{FFR}.
In \cite[Sect.~11]{cris}, it is proven that ${\mathcal C}=\lim\limits_{u\to 0}\Ac_{\mu(u)}$ has 
 a simple spectrum in any
irreducible finite-dimensional $\gt g$-module $V$.  

\begin{ex}
Suppose that $\gt g=\gt{sp}_{2n}$.
The subalgebras $\widetilde{\mathcal C}$ admit an explicit description.
For a particular choice of the parameters $h(i)$, we obtain a symplectic analogue ${\mathcal{GT}}(\spa_{2n})$ of
the Gelfand--Tsetlin subalgebra.
In the notation of Section~\ref{subsec:genbcd} take
\ben
h(i{-}1)=F_{i\tss i}\qquad\text{for}\quad i=1,\dots,n.
\een
Then $\gt g_{m{-}1}=\mK^{m}\oplus\spa_{2n{-}2m}$, where
we identify $\spa_{2n{-}2m}$ with the Lie subalgebra of $\spa_{2n}$
spanned by the elements $F_{ij}$ with $m+1\le i,j\le (m+1)'$.
The arising algebra ${\mathcal{GT}}(\spa_{2n})=\widetilde{\mathcal C}$ is freely generated by
the centres $\U(\spa_{2k})^{\spa_{2k}}$ with $1\le k\le n$ and by
the symmetrisations
\ben
\varpi\big(\partial_{h(m)} {\Phi^{(m)}_{2i}}\big)\in\U(\gt{sp}_{2n{-}2m})
\qquad
\text{for} \quad m=0,\dots,n-1 \fand i=1,\dots,n-m,
\een
where $\Phi_{2i}^{(0)}=\Phi_{2i}$ and
$\Phi_{2i}^{(m)}\in\cS(\gt{sp}_{2n-2m})\subset \cS(\gt g_{m{-}1})$ with $m\ge 1$
denotes the
$2i$-th coefficient of
the characteristic polynomial associated with $\spa_{2n-2m}$
as defined in \eqref{detsp}.
The subalgebra ${\mathcal{GT}}(\spa_{2n})$ is maximal commutative.
Making use of the anti-involution on ${\mathcal U}\big(t^{-1}\gt g[t^{-1}]\big)$
defined in the proof of \cite[Lemma~2]{FFR}, one shows that
the action of  ${\mathcal{GT}}(\spa_{2n})$ on 
any finite-dimensional irreducible $\spa_{2n}$-module 
$V$ is diagonalisable \cite{cris}.
Moreover, 
this subalgbera 
has  a simple
spectrum in $V$ \cite{cris},
thus providing  a basis of Gelfand--Tsetlin-type; cf.~\cite{m:br}.

In the first nontrivial example with $n=2$ the subalgebra ${\mathcal{GT}}(\spa_{4})\subset\U(\spa_4)$
is generated by the centres of $\U(\spa_2)$ and $\U(\spa_4)$, the Cartan elements
$F_{11},F_{22}$ and one more element
\ben
\Det\begin{bmatrix}F_{22}&F_{23}&F_{24}\\
F_{32}&F_{33}&F_{34}\\
F_{42}&F_{43}&F_{44}
\end{bmatrix}
-\Det\begin{bmatrix}F_{11}&F_{12}&F_{13}\\
F_{21}&F_{22}&F_{23}\\
F_{31}&F_{32}&F_{33}
\end{bmatrix},
\een
where we used the symmetrised determinants introduced in Section~\ref{subsec:detper}.
\end{ex}

\subsection{Vinberg's problem for centralisers}
Let $\gamma\in\gt g\cong\gt g^*$ be a nilpotent element. Set $\gt q=\gt g_\gamma$, take $\nu\in\gt q^*$ and
consider the Mishchenko--Fomenko subalgebra $\oA_\nu\subset\cS(\gt q)$.
In this setting, Vinberg's problem was recently solved by Arakawa and Premet \cite{ap}
under the assumptions that $\gt q$ has the codim-2 property, there is a g.g.s. for $\gamma$, and
$\nu\in\gt q^*_{\rm reg}$. We will restrict ourselves to types {\sf A} and {\sf C}, where the first two  assumptions are satisfied.

Let $\{H_1,\ldots,H_n\}\subset\cS(\gt g)^{\gt g}$ be the set of generators, where
$H_i=\Phi_{i}$ in type {\sf A} and $H_i=\Phi_{2i}$ in type {\sf C}.
Recall that $H_1,\ldots,H_n$ is a g.g.s. for any nilpotent $\gamma\in\gt g$ \cite{ppy}.
Set $P_i:={^\gamma\!H_i}$. Then $\cS(\gt q)^{\gt q}=\mK[P_1,\ldots,P_n]$ by  \cite{ppy}.

\begin{prop}\label{c-kom}
Suppose that $\gt g$ is of type {\sf A} or {\sf C}.
Take any $\nu\in\gt q^*$. Then the elements $\varpi(\partial_\nu^k P_i)$ generate a commutative
subalgebra $\widetilde{\Ac}_\nu$ of $\U(\gt q)$. If $\nu\in\gt q^*_{\rm reg}$, then
$\gr\widetilde{\Ac}_\nu=\oA_\nu$.
\end{prop}
\begin{proof}
Suppose first that $\nu\in\gt q^*_{\rm reg}$ and that $\nu=\bar\mu=\mu|_{\gt q}$ with $\mu\in\gt g^*_{\rm reg}$.  Then
$$\varpi(\partial_{\gamma{+}u\mu}^k H_i)\in {\mathcal A}_{\gamma{+}u\mu}$$ as in the proof of
Proposition~\ref{cor:limvin}.
Recall that if $k\ge \deg P_i$ and $\bar k=k-\deg P_i$, then
$$\lim\limits_{u\to 0}\mK \partial_{\gamma{+}u\mu}^kH_i = \mK \partial_{\nu}^{\bar k}P_i.$$
 Hence
$\varpi(\partial_\nu^k P_i)\in{\mathcal C}_{\gamma,\mu}$ for all $k\ge 0$  and therefore any two such elements commute.
The statement holds for all $\nu$ in a dense open subset
implying that
\ben
\big[\varpi(\partial_{\nu'}^k P_i),\varpi(\partial_{\nu'}^{k'} P_{i'})\big]=0
\een
for all $\nu'\in\gt q^*$, all $k,k'\in\mathbb Z_{\ge 0}$ and all
$i,i'\in\{1,\ldots,n\}$.

The Lie algebra $\gt q$ satisfies  the codim-2 condition and has  $n$, where  $n=\ind\gt q$, algebraically independent
symmetric invariants  \cite{ppy}. Hence ${\rm tr.deg}\,\oA_\nu={\bf b}(\gt q)$ if $\nu\in\gt q^*_{\rm reg}$ \cite[Thm~3.1]{bol}, see also \cite[Sect.~2.3]{codim3}. At the same time
 $\sum\limits_{i=1}^n \deg P_i={\bf b}(\gt q)$ \cite{ppy}.
Therefore the elements $\partial^{\tss k}_{\nu} P_i$ with $0\le k<\deg P_i$ have to be algebraically independent
for each $\nu\in\gt q^*_{\rm reg}$; see also Proposition~\ref{Prop-Bol}. 
By a standard argument, the symbol of ${\bf Q}\big(\varpi(\partial_\nu^k P_i)\big)$ is equal to
${\bf Q}(\partial_\nu^k P_i)$ for every polynomial ${\bf Q}$ in ${\bf b}(\gt q)$ variables.
Hence $\gr\widetilde{\Ac}_\nu=\oA_\nu$.
\end{proof}

\begin{conj}\label{conj:ap}
If $\nu\in\gt q^*_{\rm reg}$, then $\widetilde{\Ac}_\nu$ coincides with the quantisation
of  $\oA_\nu$ constructed in \cite{ap}.
\end{conj}

For a reductive $\gt g$, the uniqueness of the quantisation
in the case of a generic semisimple  element $\mu\in\gt g\cong \gt g^*$ is proven by Rybnikov \cite{L2}.
However, it is not known whether this uniqueness property extends to the quantisation
of $\oA_\nu$ and we cannot conclude that  the symmetrisation in the sense of
Proposition~\ref{c-kom} coincides with the quantisation of Arakawa and Premet.

\end{document}